\documentclass[12pt]{article}

\usepackage{lmodern}
\usepackage[T1]{fontenc}
\usepackage[english]{babel}

\usepackage[a4paper, top=2cm,bottom=2cm,left=3cm,right=3cm,marginparwidth=1.75cm]{geometry}

\usepackage{mathtools} 
\usepackage{amssymb}
\usepackage{amsfonts}
\usepackage{amsthm}
\usepackage{nicefrac,xfrac}

\usepackage{subcaption} 
\usepackage{enumitem}
\usepackage{booktabs}

\usepackage[hyphens]{url}
\usepackage[colorlinks=true, allcolors=blue]{hyperref}

\graphicspath{{./figs/}}

\newtheorem{theorem}{Theorem}[section]
\newtheorem{lemma}[theorem]{Lemma}
\newtheorem{corollary}[theorem]{Corollary}
\newtheorem{proposition}[theorem]{Proposition}

\newtheorem{assumption}{Assumption}[section]
\newtheorem{remark}{Remark}[section]

\DeclarePairedDelimiter{\norm}{\lVert}{\rVert}

\newcommand{\E}{{\mathbb{E}}}
\newcommand{\EM}{\mathop{\mathbb{E}_\M}}

\newcommand{\R}{\mathbb{R}}
\newcommand{\M}{\mathcal{M}}

\DeclareMathOperator{\br}{Br}
\DeclareMathOperator{\ibr}{IBr}
\DeclareMathOperator{\mise}{MISE}
\DeclareMathOperator{\mse}{MSE}
\newcommand{\ibias}{\operatorname{IBias}^2}
\newcommand{\bias}{\operatorname{Bias}^2}
\DeclareMathOperator{\ivar}{IVar}
\DeclareMathOperator{\var}{Var}
\DeclareMathOperator{\Rot}{Rot}

\begin{document}

\title{Bias Reduction by Multiscale Quasi-Interpolation for Scalar- and Manifold-Valued Functions}

\author{
Asaf Abas \quad and \quad Nir Sharon \\
Department of Applied Mathematics\\
Tel Aviv University\\
Tel Aviv-Yafo, Israel
}
\date{}

\maketitle

\begin{abstract}
We study the bias--variance tradeoff within a multiscale approximation framework. Our approach uses a given quasi-interpolation operator, which is repeatedly applied within an error-correction scheme over a hierarchical data structure. We use the bias ratio, the fraction of mean squared error attributable to squared bias, to assess multiscale estimators. We also provide a theoretical mechanism for the observed bias reduction: for scalar-valued linear quasi-interpolation with zero-mean additive noise, the statistical bias of the multiscale estimator is exactly the deterministic residual produced by the multiscale error-correction scheme. A residual-contraction condition then yields bias reduction, and a spectral-filter model gives a monotone nonincreasing integrated bias ratio with respect to the number of levels. Within the same model, we obtain an exact per-level error increment whose bias--variance crossover gives a criterion for when an additional level is beneficial. Under a local consistency assumption, a tangent-space analysis extends the residual-bias mechanism to manifold-valued constructions and yields a leading-order bias--variance decomposition. Our findings establish multiscale approximation as a bias-reduction methodology applicable to general quasi-interpolation operators, including applications to manifold-valued functions.
\end{abstract}


\section{Introduction}

Quasi-interpolation is a standard approach for the smooth approximation of functions given as samples over scattered data sites~\cite{wendland2004scattered, levin1998approximation}. This type of operator has been generalized for various domains, particularly for functions defined on manifolds, e.g.,~\cite{grohs2017scattered, grohs2020handbook, zimmermann2021manifold, sharon2022multiscale}. Quasi-interpolation methods are known for reproducing specific function spaces, preserving stability and locality, and maintaining low computational complexity. However, the smooth results they produce often yield high-bias approximations.

To further improve these properties, we use the standard error-correction multilevel (multiscale) scheme as a black-box wrapper around a given quasi-interpolation operator to mitigate high-bias approximations. This strategy uses nested sets from the data and applies error correction across multiple scales, enabling more effective capture of finer details. The multilevel error-correction approach was initially introduced in the context of radial basis functions~\cite{floater1996multistep}; this hierarchical error-correction scheme can be readily extended to many types of approximation operators.

Recent studies on multiscale approximation have focused on analyzing the error rates associated with the multiscale approach~\cite{wendland2010multiscale, franz2023multilevel}, while others have adapted it to new settings~\cite{sharon2022multiscale} or made improvements~\cite{le2020multilevel}. Most recently, Durst and Wendland~\cite{durst2026improved} proved that the multilevel moving least squares (MLS) scheme of~\cite{sharon2022multiscale} converges, when the data sites form a regular grid, at a rate that exceeds the rate of single-scale MLS of the same reproduction degree~$r$ by up to nearly~$r+1$ additional orders for sufficiently smooth targets, thereby giving a rigorous account of the accelerated convergence observed numerically in~\cite{sharon2022multiscale}. A line of research has expanded the Nystr\"{o}m extension with multiscale error correction schemes~\cite{bermanis2013multiscale,duchateau2013adaptation}. This approach has applications in dimension reduction using diffusion maps~\cite{rabin2012heterogeneous,chiavazzo2014reduced}. In contrast to our setting, the Nystr\"{o}m extensions exclusively utilize the Gaussian kernel. Most of this literature emphasizes deterministic approximation error rates or extensions to new settings; here we focus instead on the bias--variance perspective under noisy sampling, and on practical diagnostics that isolate bias reduction.

A closely related statistical perspective is provided by iterative bias reduction and $L_2$-boosting \cite{cornillon2017iterative}. For a fixed linear smoother $S$, repeated residual correction produces the estimator $[I-(I-S)^n]Y$, with the number of iterations controlling the bias--variance tradeoff. The present work examines this mechanism for multiscale quasi-interpolation with level-dependent operators on nested scattered data sets, and for manifold-valued constructions. Our contribution is the combined analysis of the ordered residual product, the integrated bias ratio and per-level error increment in a common spectral model, and the local manifold extension under an explicit consistency assumption. The scalar residual-bias identity connects these settings to the established bias-correction interpretation.

This paper studies the bias--variance behavior of multiscale quasi-interpolation for scalar- and manifold-valued functions. For linear quasi-interpolation, we identify the conditional bias with the deterministic multiscale residual, allowing deterministic approximation bounds to transfer directly to bias bounds. In a common spectral model, we establish monotonicity of the integrated bias ratio and derive an exact per-level error increment that characterizes when an additional level reduces MISE. Under an explicit local consistency assumption, we extend the residual-bias interpretation to manifold-valued constructions and obtain a leading-order bias--variance decomposition. Numerical experiments with MLS, manifold-valued Shepard approximation, and diffusion filters illustrate these results, distinguishing empirical behavior for scattered-data schemes from direct verification of the spectral model.

\section{The multiscale approximation and the bias--variance tradeoff}\label{sec:setting}

We focus on approximating scattered data. Specifically, let $f \colon \R^d \to \R$ be our target function to approximate. The data consists of the set of tuples $X = \{ (x_i, f_i) \}_{i=1}^N$ where $\{ x_i \}_{i=1}^N  \subseteq \mathbb{R}^d$ are the parametric data sites, and $\{ f_i \}_{i=1}^N$ are the observable samples. The data sites are not limited to a grid-like formation or by any other restriction. Noise may contaminate these observations, $f_i \approx f(x_i)$. The goal is to approximate the value of the function~$f(x^\ast)$ at a new data site~$x^\ast \in \R^d$.

In the manifold setting,  we have the target function $F \colon \R^d \to \M$ where $\M$ is a known manifold. We assume that the geodesic distance~$\rho:\M \times \M \to \R$ is available in the sense that the logarithm map is well-defined on~$\M$ (and forms a unique shortest geodesic) for the point pairs encountered in our computations. This holds globally on Hadamard manifolds (also known as Cartan--Hadamard manifolds) as simply connected, complete, and non-positive sectional curvature manifolds, and the geodesic between any two points is unique~\cite{li2012geometric}. Nevertheless, it is also valid locally on general manifolds, e.g., using the principal logarithm on points close enough, as done, for example, on SO(3).

Another limitation in the manifold setting is the need to approximate the expected value numerically using the Karcher mean~\cite{karcher1977riemannian}. For Hadamard manifolds, the Karcher mean exists and is unique~\cite{huning2019convergence, sander2016geodesic}. For manifolds with positive curvature, such as in our test case of SO(3), we can pose conditions on the spreading of points; see~\cite{dyer2016barycentric, huning2022convergence}. Regarding our measurements, which are possibly noisy, we introduce noise in the tangent space of the manifold. We compute the noisy sample as~$F_i = \exp_{F(x_i)}(\xi_i)$, where~$\xi_i$ is the noise term, sampled from a Gaussian in the tangent space, and~$\exp$ is the Riemannian exponential map on~$\M$.

Next, we detail the multiscale framework for both scalar-valued and manifold-valued functions. Additionally, we present the bias--variance tradeoff and propose a new metric based on it, the bias ratio.

\subsection{The multiscale approximation}

Our multiscale (or multilevel) method serves as an error-correction scheme and was first introduced in~\cite{floater1996multistep}. This method can enhance the results of a given operator without any additional data or assumptions~\cite{sharon2022multiscale, wendland2010multiscale, franz2023multilevel, le2020multilevel, bermanis2013multiscale, chiavazzo2014reduced}. The multiscale approximation is based on two main procedures: generating a sequence of nested hierarchical sets and applying a given operator~$Q$ via error-correction schemes.

A hierarchical sequence of subsets of~$ X$ satisfies~$ X_1 \subset \ldots \subset X_n = X$. Each subset~$X_i$ corresponds to a scale (or level) of the multiscale. The number of scales~$n$ is a parameter. Throughout this work, we generate the subsets randomly. Their sizes satisfy a given proportional growth rate~$\gamma := {\#X_{i-1}}  /{\#X_i}$, so $\gamma \in (0,1)$. Additionally, the multiscale scheme involves deploying a given ``black-box'' approximation operator~$Q$. We focus on quasi-interpolation operators. We denote the approximation  of~$Q$ over the dataset~$X_j$ with the target function~$f$ by~$Q_j f = Q_{X_j} f$.

The error correction  scheme defines a sequence of operators~$\{M_i\}_{i=1}^n$, each associated with a corresponding subset~$\{X_i\}_{i=1}^n$. The literature presents two equivalent approaches for expressing the multiscale operators~$M_i$: an iterative approach~\cite{floater1996multistep} and a direct approach~\cite{franz2023multilevel}. While the direct approach is commonly used for theoretical analysis, the iterative approach offers practical advantages for numerical implementation. The iterative approach constructs two sequences of operators: the multiscale approximation operators~$\{M_i\}_{i=1}^n$ and the multiscale error operators~$\{E_i\}_{i=1}^n$. The initial conditions are~$M_0 = 0$ and~$E_0 f = f$. At each step~$i$, the method updates the operators according to:
\begin{equation} \label{eq:multiscale-recursion}
\begin{split}
    M_i f \, = \, M_{i-1} f + Q_iE_{i-1}f, \\
    E_i f \,= \,\,  E_{i-1}f - Q_i E_{i-1}f.
    \end{split}
\end{equation}
After completing~$n$ steps, this process produces the approximated function~$\widetilde{f} = M_n f$.

The numerical implementation of the multiscale for scalar-valued functions is straightforward~\eqref{eq:multiscale-recursion}. In~\cite{sharon2022multiscale}, the authors discuss the generalizations that are needed to approximate a manifold-valued function, and we adopt their modifications in our implementation. These modifications include replacing addition and subtraction operations with the exponential map and logarithmic map, respectively.

\subsection{Bias--variance tradeoff for scalar- and manifold-valued functions}\label{sec:bvt}

The bias--variance tradeoff (BVT) is extensively studied and applied across various fields, such as statistics, signal processing, and more. Notable applications of the BVT include statistical analysis~\cite{kazashi2021density, kohavi1996bias} and algorithm tuning~\cite{geman1992neural, zhang2022mitigating}. This section introduces the BVT in the context of approximating scalar-valued functions and manifold-valued functions.

Let~$\widetilde{f}$ be an approximation of a scalar-valued target function~$f \colon \R^d \to \R$. The mean squared error ($\mse$) of~$\widetilde{f}$ is
\begin{equation} \label{eq:mse}
    \mse[\widetilde{f},f](x) = \E\left[(f(x) - \widetilde{f}(x))^2\right],
\end{equation}
where~$\E\left[\cdot\right]$ denotes the expected value with respect to the sampling of the domain and the noise. We will need to distinguish this unconditional expectation from the conditional one, taken with respect to the noise only for a fixed set of data sites and a fixed hierarchy. The theory of Section~\ref{sec:linear} is developed conditionally, which is the sharper statement; Proposition~\ref{prop:random-design} then transfers it to the unconditional quantities~\eqref{eq:mse} and~\eqref{eq:bias-var} that the numerical protocol of Section~\ref{sec:num-settings} estimates, and identifies the extra, design-induced variance component that appears there. The BVT decomposes the MSE into two components: bias and variance. The MSE decomposes as
\begin{equation} \label{eq:bvt}
    \mse[\widetilde{f},f](x)  = \bias[\widetilde{f},f](x) + \var[\widetilde{f}](x) ,
\end{equation}
where
\begin{equation}\label{eq:bias-var}
\begin{split} 
    \bias[\widetilde{f},f](x)&=\left(\E \left[ f(x) - \widetilde{f}(x)\right]\right)^2= \left( f(x)  - \E \left[ \widetilde{f}(x) \right] \right)^2, \\
    &\var[\widetilde{f}](x) = \E  \left[\left(\widetilde{f}(x) - \E [\widetilde{f}(x)]\right)^2\right].
\end{split} 
\end{equation}
As in most scenarios, we cannot compute these terms directly and must approximate them. To this end, we perform multiple trials. From~$T$ different trials, we compute a set of approximations~$\{ \widetilde{f}^{1},\ldots,  \widetilde{f}^{T}\}$. We use their average as an estimate of the expected value,
\begin{equation} \label{eq:mc-mean}
    \E \left[\widetilde{f}(x)\right] \approx \frac{1}{T} \sum_{j=1}^{T} \widetilde{f}^j(x).
\end{equation}

Next, we revise the above definitions for the case of a manifold-valued function, $F \colon \R^d \to \M$. The extension of~\eqref{eq:mse} and~\eqref{eq:bias-var} is achieved using the geodesic distance~$\rho$ associated with the manifold~$\M$. In particular, it is common practice when computing the MSE on a manifold to measure the error using the squared geodesic distance, e.g.,~\cite{jermyn2005invariant}. Therefore, we define the MSE of an approximation~$\widetilde{F}$ of the manifold-valued target function~$F$ as
\begin{equation}\label{eq:mse-manifold}
   \mse[\widetilde{F},F](x) = \E \left[ \rho(F(x), \widetilde{F}(x))^2 \right].
\end{equation}
For completeness, we present the remaining definitions:
\begin{equation} \label{eq:bias-manifold}
\bias[\widetilde{F},F](x) = \rho\left(F(x), \EM[\widetilde{F}(x)]\right) ^2,
\end{equation}
\begin{equation} \label{eq:var-manifold}
\var[\widetilde{F}](x) = \E\left[\rho\left(\EM[\widetilde{F}(x)], \widetilde{F}(x)\right)^2\right].
\end{equation}
To extend~\eqref{eq:mc-mean} to manifold-valued functions, we introduce~$\EM$ to denote the expected value on the manifold~$\M$. An exact global bias--variance decomposition does not hold on manifolds; however, a leading-order decomposition is available (Lemma~\ref{lem:manifold-variance}), and beyond that leading order, the bias ratio is used as a diagnostic tool.

\subsection{The bias ratio}\label{sec:bias-ratio}

While the multiscale scheme itself is standard, the goal here is to quantify its effect on the bias component of the error via a normalized diagnostic. We focus on the effect of multiscale approximation on the bias component in the MSE. In general, multiscale error correction behaves like residual refinement, which increases the effective approximation order for smooth functions (see~\cite{durst2026improved} for a rigorous statement for MLS on regular grids); hence, the bias decreases faster. To enable comparison of the bias across different approximation operators without being affected by their varying MSE values, we suggest the bias ratio:
\begin{equation}\label{eq:bias-ratio}
    \br[\widetilde{f},f](x) = \frac{\bias[\widetilde{f},f](x)} {\mse[\widetilde{f},f](x)}.
\end{equation}
In the Euclidean setting, where~\eqref{eq:bvt} holds, we have that~$\br\in[0,1]$, with~$\br=1$ if and only if~$\var=0$ (assuming $\mse>0$). For manifold-valued outputs, we adopt the analogs~\eqref{eq:bias-manifold} and~\eqref{eq:var-manifold} of bias and variance to obtain $\br[\widetilde{F}, F](x)$ for a manifold-valued function $F$. In this case an exact decomposition does not hold; instead, we use $\br$ as a proxy for the ratio of the squared bias in the MSE, justified to leading order by Lemma~\ref{lem:manifold-variance}. Note that, in our experiments, we empirically observe that $\br$ takes values in $[0,1]$.

We stress what this diagnostic is, and is not, intended for. The bias ratio is a normalized quantity for comparing the bias content of different approximation operators, or of one operator across noise levels and dataset sizes, without being confounded by the different magnitudes of their MSE. It is deliberately not a criterion for choosing the number of levels: within the spectral-filter model of Section~\ref{sec:spectral}, the integrated bias ratio is monotone nonincreasing in the number of levels (Theorem~\ref{thm:spectral-damping}), so it can never signal a level at which to stop, and level selection is an MSE consideration governed by Corollary~\ref{cor:mise-increment}; see Remark~\ref{rem:optimal-levels}. For the same reason, a decrease in~$\br$ certifies a reduction of the bias only when it is read together with the MSE, since $\br$ also decreases when the variance grows. We therefore report both quantities throughout Section~\ref{sec:numerics}.

\section{A residual-bias interpretation of the multiscale scheme}\label{sec:residual-bias}

The notation used in this section and in the remainder of the paper is collected in Table~\ref{tab:notation}.
\begin{table}[!ht]
\centering\small
\begin{tabular}{@{}lp{0.7\textwidth}@{}}
\toprule
Symbol & Meaning \\
\midrule
$X=\{x_i\}_{i=1}^N$ & data sites; $X_1\subset\cdots\subset X_n=X$ is the hierarchy and $N_i=\#X_i$ \\
$\mathbf X$ & the design: the sites, the hierarchy, and the level parameters they induce \\
$\gamma$ & proportional growth rate of the hierarchy \\
$\Gamma$ & external evaluation grid \\
$\mathcal H=\R^{X}$ & space of scalar values on the data sites \\
$Q_i$, $R_i=I-Q_i$ & level-$i$ quasi-interpolation operator and its residual operator \\
$M_j$, $E_j$, $\mathcal R_j=R_j\cdots R_1$ & multiscale approximation, error, and accumulated residual operators \\
$h_i$, $\delta_i=\nu h_i$, $K$ & fill distance, support radius and radial kernel \\
$\varepsilon$, $\sigma^2$ & scalar observational noise and its variance \\
$\eta_i$ & per-level residual contraction factor \\
$L=U\Lambda U^{*}$, $\omega_k$, $u_k$ & operator of the spectral model, its eigenvalues and eigenvectors \\
$\phi_i$, $\lambda_{i,k}$, $a_{j,k}$ & spectral filter, spectral response, accumulated residual factor \\
$\theta_k$ & spectral coefficient of the target \\
$\M$, $\rho$, $q=\dim\M$ & manifold, geodesic distance, and dimension \\
$\xi_i$ & tangent-space noise on the manifold samples \\
$p$, $u=\log_pF$, $r$ & chart centre, normal coordinates, neighbourhood radius \\
$A_n$, $A_n(x)$ & evaluation operator on $\Gamma$ and its tangent-space analogue \\
$e_n$, $g_n$, $\chi=r^3+\sigma^2$ & evaluation residual (scalar in Section~\ref{sec:linear}, tangent-space in Section~\ref{sec:manifold}), consistency remainder, remainder scale \\
$\zeta$, $\Sigma$, $\beta$ & centred coordinate noise, its covariance, and its mean shift \\
$\mathcal B_n^2$, $\mathcal V_n$ & integrated tangent-space squared bias and variance \\
\bottomrule
\end{tabular}
\caption{Notation used throughout the paper. Symbols not listed here are local to the statement or proof in which they appear.}
\label{tab:notation}
\end{table}

We now give a theoretical explanation for why the multiscale construction in~\eqref{eq:multiscale-recursion} acts as a bias-reducing mechanism. It is important to distinguish three tiers of results, which rest on progressively stronger assumptions and which cover the experiments to different degrees. The first tier assumes \emph{linearity} of the multiscale operators. Lemma~\ref{lem:direct-residual-form} and Theorem~\ref{thm:bias-equals-residual} show that the statistical bias of the multiscale estimator equals the deterministic multiscale residual~$\mathcal R_nf$, and Corollary~\ref{cor:residual-contraction} converts any residual contraction into bias reduction. This tier applies verbatim to fixed-site moving least squares (MLS) of any polynomial degree, since such operators are exactly linear in the data; no self-adjointness, common eigenbasis, or spectral constraint is required.

The second tier is a \emph{spectral-filter model} (Theorem~\ref{thm:spectral-damping} and Corollary~\ref{cor:mise-increment}) in which the level operators additionally share an orthonormal eigenbasis with responses in~$[0,1]$. This buys the extra conclusions of monotone bias-ratio decay and a bias--variance crossover, at the cost of assumptions that general scattered-data operators need not satisfy exactly.

The third tier is a \emph{local tangent-space result} for manifold-valued functions (Lemmas~\ref{lem:manifold-linearization} and~\ref{lem:manifold-variance}), which transfers the first-tier mechanism to the nonlinear manifold schemes to leading order. Throughout, the assumptions are sufficient rather than necessary; their purpose is to isolate the operator-level mechanism behind the bias-ratio reductions observed in the experiments.

As we make precise below, the exact backbone (first tier) already covers the scalar MLS experiments: Corollary~\ref{cor:residual-contraction} turns any per-level contraction of the deterministic residual into a bias bound. For multilevel MLS on gridded data, such contraction factors are available from the recent analysis of Durst and Wendland~\cite{durst2026improved} (Remark~\ref{rem:contraction-hypothesis}); for the random-site, Gaussian-weighted MLS used in our experiments, they remain a hypothesis that the experiments support indirectly through the observed bias reduction. The second tier is realized exactly by a class of diffusion-type operators, for which Section~\ref{sec:diffusion} provides a numerical experiment.

\subsection{Multiscale with linear quasi-interpolation operators}\label{sec:linear}

Throughout this section we condition on fixed data sites and a fixed hierarchy~$X_1\subset\cdots\subset X_n=X$, so that the level operators~$Q_i=Q_{X_i}$ are deterministic and~$\E$ denotes expectation with respect to the observational noise only. We take as state space the finite-dimensional Hilbert space~$\mathcal H=\R^{X}$ of scalar values on the full data set~$X=X_n$, with the Euclidean inner product~$\langle\cdot,\cdot\rangle$ and norm~$\|\cdot\|$ (a quadrature-weighted inner product can be used instead, in which case the covariance assumption below is to be read as an operator identity on~$\mathcal H$). Each level operator is the endomorphism~$Q_i:\mathcal H\to\mathcal H$ that restricts a function's values to the sub-site set~$X_i$, forms the level-$i$ quasi-interpolant, and re-evaluates it on all of~$X$; the observational noise is a random element of~$\mathcal H$.

With this convention~$Q_i$,~$M_n$, and~$\mathcal R_n$ all act on the single space~$\mathcal H$, so that the accumulation identity~$M_n=I-\mathcal R_n$ is well posed and the adjoint below is taken with respect to~$\langle\cdot,\cdot\rangle$.

Reporting the error at points outside~$X$ (for example, on a fine evaluation grid) amounts to a fixed, deterministic evaluation of the resulting approximant; since that evaluation is applied after the scheme and is independent of the noise, it commutes with the expectation.

Similar to~\cite{franz2023multilevel}, for each level,~$Q_i$ is a deterministic linear quasi-interpolation operator, and we write
\[
    R_i = I-Q_i,
    \qquad
    \mathcal R_0=I,
    \qquad
    \mathcal R_j = R_j R_{j-1}\cdots R_1,\quad j=1,\ldots,n,
\]
where~$I$ is the identity operator. Thus,~$R_i$ is the residual operator at level~$i$, while~$\mathcal R_j$ is the residual operator accumulated through level~$j$.

\begin{lemma}[Direct residual form]\label{lem:direct-residual-form}
Let~$Q_1,\ldots,Q_n$ be linear quasi-interpolation operators. Then, for~$j=0,\ldots,n$, the multiscale error and approximation operators defined in~\eqref{eq:multiscale-recursion} satisfy
\begin{equation}\label{eq:direct-error-form}
    E_j = \mathcal R_j,
\end{equation}
and
\begin{equation}\label{eq:direct-approx-form}
    M_j = I-\mathcal R_j.
\end{equation}
In particular,~$E_n= \mathcal R_n=(I-Q_n)(I-Q_{n-1})\cdots(I-Q_1)$.
\end{lemma}

\begin{proof}
The recursion in~\eqref{eq:multiscale-recursion} gives
\[
    E_i f = E_{i-1}f-Q_iE_{i-1}f=(I-Q_i)E_{i-1}f = R_iE_{i-1}f.
\]
Since~$E_0=I=\mathcal R_0$, induction gives~$E_j=R_jR_{j-1}\cdots R_1=\mathcal R_j$ for all~$j$. This is the direct form of the multiscale error operator used in~\cite{franz2023multilevel}.

Summing the expressions for~$M_i$ and~$E_i$ in~\eqref{eq:multiscale-recursion} gives
\[
    M_i f+E_i f
    =M_{i-1}f+Q_iE_{i-1}f+E_{i-1}f-Q_iE_{i-1}f
    =M_{i-1}f+E_{i-1}f.
\]
Because~$M_0=0$ and~$E_0=I$, another induction yields~$M_i+E_i=I$ for all~$i$. Hence~$M_j=I-E_j=I-\mathcal R_j$.
\end{proof}

The next theorem connects the deterministic residual in Lemma~\ref{lem:direct-residual-form} to the statistical bias and variance in~\eqref{eq:bias-var}. It is the main theoretical anchor for the paper's scalar-valued part. In the theorem, we also use integrated versions of the squared bias, variance, and MSE, denoted~$\ibias$,~$\ivar$, and~$\mise$, which are obtained by summing the pointwise quantities over the data sites with respect to the inner product of~$\mathcal H$. Throughout, ``integrated'' refers to this sum; we keep the customary acronym $\mise$ although no normalization by the number of sites is applied, so the integrated quantities scale with~$\#X$.

\begin{theorem}[Bias equals the deterministic multiscale residual on the data sites]\label{thm:bias-equals-residual}
All quantities below are taken on the data-site space~$\mathcal H=\R^{X}$ introduced at the beginning of Section~\ref{sec:linear}. Let~$Y=f+\varepsilon$ be the noisy observation of~$f\in\mathcal H$, where~$\E[\varepsilon]=0$. Let~$\widetilde f_n=M_nY$ be the scalar-valued multiscale estimator after~$n$ levels, with the hierarchy fixed as above. Then,
\begin{equation}\label{eq:expected-multiscale-estimator}
    \E[\widetilde f_n]=M_nf,
\end{equation}
and the conditional bias is
\begin{equation}\label{eq:bias-residual-identity}
    f-\E[\widetilde f_n]=E_nf=\mathcal R_nf.
\end{equation}
Consequently, at any data site~$x \in X$, the conditional squared bias is
\begin{equation}\label{eq:pointwise-bias-residual}
   \left(f(x)-\E[\widetilde f_n(x)]\right)^2=\left((\mathcal R_nf)(x)\right)^2.
\end{equation}

If, in addition,~$\operatorname{Cov}(\varepsilon)=\sigma^2 I$, then the conditional integrated quantities satisfy
\begin{equation}\label{eq:integrated-bias--variance-residual}
\begin{split}  
    \ibias_n &= \|\mathcal R_nf\|^2,
    \quad
    \ivar_n = \sigma^2\operatorname{tr}(M_nM_n^*),
     \\
    \mise_n &= \|\mathcal R_nf\|^2+\sigma^2\operatorname{tr}(M_nM_n^*),
    \end{split}
\end{equation}
where~$M_n^*$ is the adjoint of~$M_n$.
\end{theorem}

\begin{proof}
By linearity  of~$M_n$ and the zero-mean noise assumption,
\[
    \E[\widetilde f_n]
    =\E[M_n(f+\varepsilon)]
    =M_nf+M_n\E[\varepsilon]
    =M_nf.
\]
Using Lemma~\ref{lem:direct-residual-form},
\[
    f-\E[\widetilde f_n]
    =f-M_nf
    =(I-M_n)f
    =E_nf
    =\mathcal R_nf,
\]
which proves~\eqref{eq:bias-residual-identity}. The pointwise identity~\eqref{eq:pointwise-bias-residual} is its squared evaluation at~$x$. For the integrated formulas, take the squared norm of~\eqref{eq:bias-residual-identity}. The variance term is
\[
    \E\left[\|\widetilde f_n-\E[\widetilde f_n]\|^2\right]
    =\E\left[\|M_n\varepsilon\|^2\right]
    =\operatorname{tr}\left(M_n\E[\varepsilon\varepsilon^*]M_n^*\right)
    =\sigma^2\operatorname{tr}(M_nM_n^*).
\]
Adding the integrated squared bias and variance gives~\eqref{eq:integrated-bias--variance-residual}.
\end{proof}
Theorem~\ref{thm:bias-equals-residual} is a statement on the data-site space~$\mathcal H=\R^{X}$: the residual~$\mathcal R_nf$ and the operator~$M_n$ are square endomorphisms of~$\mathcal H$. In the experiments, however, the error is reported on an external evaluation grid, by evaluating the multiscale approximant there and comparing it to the target. This uses a different, rectangular operator, which we treat next.

We fix a finite evaluation grid~$\Gamma$ and let~$A_n:\mathcal H\to\R^{\Gamma}$ be the deterministic linear map that sends the noisy data~$Y = f + \varepsilon$ to the multiscale approximant evaluated on~$\Gamma$. The map~$A_n$ shares the level fits of~$M_n$ but is rectangular when~$\Gamma\neq X$.

\begin{proposition}[Evaluation on an external grid]\label{prop:grid-evaluation}
Denote by~$f_X=f|_X\in\mathcal H$ and~$f_\Gamma=f|_G\in\R^{\Gamma}$ the restrictions of the target function, and assume~$\E[\varepsilon]=0$ with $\operatorname{Cov}(\varepsilon)=\sigma^2I$ as in Theorem~\ref{thm:bias-equals-residual}. Then, the estimator~$\widetilde f_n^{\,\Gamma}=A_nY$ satisfies
\[
    \E[\widetilde f_n^{\,\Gamma}]=A_nf_X ,
\]
and, for every~$x\in\Gamma$,
\[
    \bias[\widetilde f_n^{\,\Gamma},f_\Gamma](x)=\big(f_\Gamma(x)-(A_nf_X)(x)\big)^2,
    \qquad
    \var[\widetilde f_n^{\,\Gamma}](x)=\sigma^2\|a_n(x)\|_2^2 ,
\]
where~$a_n(x)\in\R^{X}$ denotes the row of~$A_n$ associated with~$x$. The corresponding integrated quantities are
\[
    \ibias[\widetilde f_n^{\,\Gamma},f_\Gamma]=\|f_\Gamma-A_nf_X\|^2,
    \qquad
    \ivar[\widetilde f_n^{\,\Gamma}]=\sigma^2\operatorname{tr}(A_nA_n^{*})=\sigma^2\|A_n\|_F^2 .
\]
Thus, the bias reported on the grid is exactly the deterministic approximation error of the noiseless scheme evaluated on~$\Gamma$. Specifically, when~$\Gamma=X$, we have~$A_n=M_n$, which reduces to Theorem~\ref{thm:bias-equals-residual}.
\end{proposition}

\begin{proof}
Since~$A_n$ is linear and deterministic and~$\E[\varepsilon]=0$, we have $\E[A_nY]=A_n\E[Y]=A_nf_X$, whence the pointwise bias is~$\big(f_\Gamma(x)-\E[(A_nY)(x)]\big)^2=\big(f_\Gamma(x)-(A_nf_X)(x)\big)^2$ and the integrated bias is~$\|f_\Gamma-A_nf_X\|^2$. For the variance, $(A_nY)(x)-\E[(A_nY)(x)]=\langle a_n(x),\varepsilon\rangle$, so that $\var[\widetilde f_n^{\,\Gamma}](x)=\E\big[\langle a_n(x),\varepsilon\rangle^2\big]=\sigma^2\|a_n(x)\|_2^2$. Summing over~$x\in\Gamma$ gives $\E\|A_n\varepsilon\|^2=\operatorname{tr}(A_n\E[\varepsilon\varepsilon^*]A_n^*)=\sigma^2\operatorname{tr}(A_nA_n^*)$, and the trace equals~$\|A_n\|_F^2$. Both formulas are well defined for the rectangular~$A_n$, whose adjoint maps~$\R^{\Gamma}\to\mathcal H$.
\end{proof}

Proposition~\ref{prop:grid-evaluation} shows that the statistical bias on the evaluation grid is exactly the deterministic grid-approximation error of the noiseless multiscale scheme. A transfer from data-site residual estimates to this external-grid error requires a separate sampling or extension estimate and is not automatic. In the numerical experiments, we therefore estimate the grid bias directly.

\begin{proposition}[Evaluation under a random design]\label{prop:random-design}
Assume that $\Gamma$ is deterministic and that all quantities
below have finite second moments. Theorem~\ref{thm:bias-equals-residual} and Proposition~\ref{prop:grid-evaluation} are stated for fixed data sites and a fixed hierarchy. Suppose instead that the design
\[
    \mathbf X=(X_1,\ldots,X_n)
\]
is random, that it determines the level operators (including any level parameters computed from it, such as the support radii), and that it is independent of the observational noise, with $\E[\varepsilon\mid\mathbf X]=0$ and $\operatorname{Cov}(\varepsilon\mid\mathbf X)=\sigma^2I$. Write $A_n(\mathbf X)$ for the resulting evaluation operator of Proposition~\ref{prop:grid-evaluation}, $a_n(x;\mathbf X)$ for its row at $x\in\Gamma$, $f_X$ for the restriction of $f$ to the realized sites, and
\[
    e_n(x;\mathbf X)
    =
    f_\Gamma(x)-\bigl(A_n(\mathbf X)f_X\bigr)(x)
\]
for the deterministic grid residual of the noiseless scheme on the realized design; on the data sites $e_n(\cdot\,;\mathbf X)=\mathcal R_n(\mathbf X)f$. Then, for the estimator $\widetilde f_n^{\,\Gamma}=A_n(\mathbf X)Y$ and every $x\in\Gamma$:
\begin{enumerate}[label=(\roman*)]
\item conditionally on the design, $f_\Gamma-\E[\widetilde f_n^{\,\Gamma}\mid\mathbf X]=e_n(\cdot\,;\mathbf X)$;
\item the unconditional bias is the averaged residual,
\[
    f_\Gamma-\E[\widetilde f_n^{\,\Gamma}]=\E_{\mathbf X}\bigl[e_n(\cdot\,;\mathbf X)\bigr],
    \qquad
    \bias[\widetilde f_n^{\,\Gamma},f_\Gamma](x)=\bigl(\E_{\mathbf X}e_n(x;\mathbf X)\bigr)^2 ,
\]
and, by Jensen's inequality, $\|f_\Gamma-\E[\widetilde f_n^{\,\Gamma}]\|\leq\E_{\mathbf X}\|e_n(\cdot\,;\mathbf X)\|$;
\item the variance splits into a noise part and a design part,
\[
    \var[\widetilde f_n^{\,\Gamma}](x)
    =
    \sigma^2\,\E_{\mathbf X}\|a_n(x;\mathbf X)\|_2^2
    +
    \var_{\mathbf X}\bigl(e_n(x;\mathbf X)\bigr);
\]
\item consequently
\[
    \mse[\widetilde f_n^{\,\Gamma},f_\Gamma](x)
    =
    \E_{\mathbf X}\bigl[e_n(x;\mathbf X)^2\bigr]
    +
    \sigma^2\,\E_{\mathbf X}\|a_n(x;\mathbf X)\|_2^2 ,
\]
so that the mean squared deterministic residual splits as
\begin{equation}\label{eq:residual-splits}
    \E_{\mathbf X}\bigl[e_n(x;\mathbf X)^2\bigr]
    =
    \bias[\widetilde f_n^{\,\Gamma},f_\Gamma](x)
    +
    \var_{\mathbf X}\bigl(e_n(x;\mathbf X)\bigr).
\end{equation}
\end{enumerate}
\end{proposition}

\begin{proof}
Part~(i) is Proposition~\ref{prop:grid-evaluation} applied conditionally on $\mathbf X$, together with $A_n(\mathbf X)=M_n(\mathbf X)$ when $\Gamma=X$. Taking expectations in~(i) and using the tower property gives the first identity in~(ii), and squaring it gives the second; Jensen's inequality applied to the norm gives the bound. For~(iii), write $\widetilde f_n^{\,\Gamma}=A_n(\mathbf X)f_X+A_n(\mathbf X)\varepsilon$ and use the law of total variance,
\[
    \var[\widetilde f_n^{\,\Gamma}](x)
    =
    \E_{\mathbf X}\Bigl[\var\bigl(\widetilde f_n^{\,\Gamma}(x)\mid\mathbf X\bigr)\Bigr]
    +
    \var_{\mathbf X}\Bigl(\E\bigl[\widetilde f_n^{\,\Gamma}(x)\mid\mathbf X\bigr]\Bigr),
\]
whose first term is $\sigma^2\E_{\mathbf X}\|a_n(x;\mathbf X)\|_2^2$ by Proposition~\ref{prop:grid-evaluation} and whose second term is $\var_{\mathbf X}(e_n(x;\mathbf X))$, since $\E[\widetilde f_n^{\,\Gamma}(x)\mid\mathbf X]=f_\Gamma(x)-e_n(x;\mathbf X)$ and $f_\Gamma(x)$ is deterministic. Part~(iv) follows by adding~(ii) and~(iii) and using $\E_{\mathbf X}[e_n^2]=(\E_{\mathbf X}e_n)^2+\var_{\mathbf X}(e_n)$, which is also~\eqref{eq:residual-splits}.
\end{proof}

Two consequences matter for reading the experiments. First, what the theory of this section controls is the deterministic residual $e_n(\cdot\,;\mathbf X)$ of each realized design. By~\eqref{eq:residual-splits}, its mean square is the reported squared bias \emph{plus} the design part of the variance, so a per-design residual contraction bounds the reported bias through~(ii), but it does not force the bias ratio to equal one.

Second, in a noiseless experiment, $\sigma=0$ and part~(iv) reduces to $\mse(x)=\E_{\mathbf X}[e_n(x;\mathbf X)^2]$: the entire mean squared error is the mean squared deterministic residual, and
\[
    \br[\widetilde f_n^{\,\Gamma},f_\Gamma](x)
    =
    \frac{\bigl(\E_{\mathbf X}e_n(x;\mathbf X)\bigr)^2}{\E_{\mathbf X}\bigl[e_n(x;\mathbf X)^2\bigr]}
    \in[0,1]
\]
is exactly the fraction of that residual which is systematic across designs rather than design-dependent. This is why the bias ratios reported for noiseless MLS in Section~\ref{sec:mls} lie strictly below one. A smaller bias ratio indicates a smaller systematic share
of the remaining mean squared residual. The situation of a fixed design, in which $\var_{\mathbf X}\equiv0$ and $\br=1$ whenever $\sigma=0$, is realized by the diffusion experiment of Section~\ref{sec:diffusion}, where all curves indeed approach $100\%$ as $\sigma\to0$.

The previous theorem identifies the statistical bias with the deterministic residual left by the multiscale approximation. Thus, any contraction estimate along the generated residuals immediately gives a bias estimate.

\begin{corollary}[Residual contraction]\label{cor:residual-contraction}
Assume the setting of Theorem~\ref{thm:bias-equals-residual}. Suppose that, for the residuals generated by the multiscale scheme,
\[
    \|E_i f\|_{\mathcal H}
    \leq
    \eta_i \|E_{i-1}f\|_{\mathcal H},
    \qquad i=1,\ldots,n.
\]
Then
\[
    \|f-\E[\widetilde f_n]\|_{\mathcal H}
    \leq
    \left(\prod_{i=1}^n \eta_i\right)\|f\|_{\mathcal H}.
\]
\end{corollary}
\begin{proof}
By Theorem~\ref{thm:bias-equals-residual},
\[
    \|f-\E[\widetilde f_n]\|_{\mathcal H}
    =
    \|E_nf\|_{\mathcal H}.
\]
Applying the assumed residual contraction recursively gives
\[
    \|E_nf\|_{\mathcal H}
    \leq
    \eta_n\|E_{n-1}f\|_{\mathcal H}
    \leq \cdots \leq
    \left(\prod_{i=1}^n\eta_i\right)\|E_0f\|_{\mathcal H}.
\]
Since~$E_0f=f$, the result follows.
\end{proof}

This corollary is only a sufficient criterion: it says that if each generated residual contracts relative to the previous one, then the multiscale bias contracts multiplicatively across levels. The condition is trajectory-wise; it need not hold as a global operator-norm bound for every quasi-interpolation operator. Its significance lies in the bridge it provides to the deterministic theory: because Theorem~\ref{thm:bias-equals-residual} identifies the bias with the \emph{deterministic} multilevel approximation error~$\mathcal R_nf=E_nf$, any convergence estimate for the multilevel quasi-interpolation error, and in particular any per-level contraction estimate, carries over directly to the statistical bias, with no spectral or self-adjointness assumptions. This covers non-normal, signed-weight operators such as MLS of arbitrary polynomial degree, without any eigenbasis or spectrum assumptions.

\begin{remark}[On the contraction hypothesis and random sampling]\label{rem:contraction-hypothesis}
The factors~$\eta_i$ are an assumption on the implemented scheme, not a consequence of it. A per-level reduction of this type is the behavior analyzed for multilevel scattered-data schemes with scaled compactly supported kernels under geometrically decreasing fill distances and sufficient target smoothness~\cite{wendland2010multiscale,franz2023multilevel}. For the multilevel MLS scheme itself, Durst and Wendland~\cite{durst2026improved} recently established such a reduction when the data sites form a regular grid (we write~$\alpha$ for the mesh ratio denoted~$\mu$ there, to avoid a clash with the measure~$\mu$ of Section~\ref{sec:manifold}): with mesh widths~$h_j=\alpha^jh_0$, $0<\alpha<1$, and support radii~$\delta_j=\nu h_j$, their recursive estimate \cite[Theorem~3.2]{durst2026improved} bounds a level-scaled Sobolev seminorm of the residual~$E_jf$ by~$C\big(\nu^{k-m}\alpha^{r+1}+\alpha^{\min\{2r+2,k\}}\big)$ times the corresponding seminorm of~$E_{j-1}f$ at the previous scale, for MLS of reproduction degree~$r$ with a compactly supported~$C^m$ kernel and a~$C^k$ target, $m>k>r$; this is a per-level contraction of the residual in the sense of Corollary~\ref{cor:residual-contraction}, in a different norm. The resulting multilevel error after~$L$ levels improves on single-scale MLS by a factor of the form~$J_{L-1}(\nu,\alpha)^{L-1}$, which becomes arbitrarily small as~$\nu\to\infty$ and~$\alpha\to0$ \cite[Theorem~3.4]{durst2026improved}. By Theorem~\ref{thm:bias-equals-residual}, these deterministic estimates are, verbatim, estimates for the bias of the multilevel MLS estimator under zero-mean noise on gridded data. Our experiments use random sites, a Gaussian weight, and random nested subsets, which are not covered by that analysis, so we treat the~$\eta_i$ as a hypothesis in this setting; the MLS experiments of Section~\ref{sec:mls} exhibit the resulting bias reduction. When the data sites and hierarchy are random and independent of the noise, conditioning on the realization~$\mathbf X$ gives $f-\E[\widetilde f_n\mid\mathbf X]=\mathcal R_n(\mathbf X)f$; averaging over~$\mathbf X$ with Jensen's inequality, as in Proposition~\ref{prop:random-design}, transfers any such bound to the averaged bias whenever the required fill-distance and quasi-uniformity conditions hold with high probability together with suitable integrability control on the complementary event.
\end{remark}

\subsection{Multiscale with a spectral filter}\label{sec:spectral}

The next result gives a componentwise interpretation of the bias-reduction mechanism within a spectral-filter model; it should be read as a model for the mechanism rather than as a structural assumption on every scattered-data quasi-interpolant. General quasi-interpolation operators, such as MLS or Shepard approximations on random data sites, need not be self-adjoint or simultaneously diagonalizable. Nevertheless, the assumptions below hold for a natural class of linear smoothing operators, namely spectral filters.

For example, given the fixed set of~$N$ data sites~$X\subset \R^d$, let~$\mathcal H=\R^{X}\cong\R^{N}$ be the space of function values on~$X$, and let~$L:\mathcal H\to\mathcal H$ be a symmetric positive-semidefinite operator, such as a graph Laplacian built from~$X$. Choose an orthonormal eigenbasis~$\{u_k\}_{k=1}^N$ of~$\mathcal H$ such that
\[
    L u_k=\omega_k u_k,
    \qquad
    0=\omega_1\leq \omega_2\leq\cdots\leq \omega_N .
\]
For each scale~$i$, define
\begin{equation}\label{eq:spectral-filter}
     Q_i= U \phi_i(\Lambda) U^{*},
    \qquad
    0\leq \phi_i(t) \leq 1 \quad\text{for all } t\in\operatorname{spec}(L).
\end{equation}
Here we use the spectral decomposition $L = U \Lambda U^{*}$, with $\Lambda=\operatorname{diag}(\omega_1,\ldots,\omega_N)$ and $U = [u_1\ u_2\ \cdots\ u_N]$. Then
\[
    Q_i u_k=\phi_i(\omega_k)u_k .
\]
A canonical example is the heat-kernel filter~$\phi_i(t)=e^{-\tau_i t}$, $\tau_i>0$. In this model, all~$Q_i$ share the same eigenvectors and have spectral responses in~$[0,1]$; if~$u_1$ is constant, they also reproduce constants since $\phi_i(0)=1$. The residual~$I-Q_i$ therefore multiplies the~$k$-th spectral component by~$1-\phi_i(\omega_k)$. The theorem below shows how this componentwise residual filtering translates into bias reduction and monotone decay of the integrated bias ratio.

\begin{theorem}[Spectral damping of bias]
\label{thm:spectral-damping}
Assume the setting of Theorem~\ref{thm:bias-equals-residual}, with \(\operatorname{Cov}(\varepsilon)=\sigma^2 I\). Suppose that \(Q_1,\ldots,Q_n\) share an orthonormal eigenbasis \(\{u_k\}_{k=1}^N\) and satisfy
\[
    Q_i u_k=\lambda_{i,k}u_k,
    \qquad
    0\leq \lambda_{i,k}\leq 1 .
\]
Let \(f=\sum_{k=1}^N\theta_k u_k\), and define
\[
    a_{j,k}:=\prod_{i=1}^j(1-\lambda_{i,k}),
    \qquad j=1,\ldots,n .
\]
Then
\[
    \ibias_j
    =
    \sum_{k=1}^N a_{j,k}^2\theta_k^2,
    \qquad
    \ivar_j
    =
    \sigma^2\sum_{k=1}^N(1-a_{j,k})^2 .
\]
Define
\[
    \ibr_j:=\frac{\ibias_j}{\ibias_j+\ivar_j}
\]
whenever the denominator is nonzero. Then~\(\ibias_j\) is nonincreasing in~\(j\), \(\ivar_j\) is nondecreasing in~\(j\), and~\(\ibr_j\) is nonincreasing in~\(j\).
\end{theorem}

\begin{proof}
Since \(Q_i u_k=\lambda_{i,k}u_k\), the residual operator \(R_i=I-Q_i\) satisfies
\[
    R_i u_k=(1-\lambda_{i,k})u_k .
\]
Therefore, for the first \(j\) multiscale levels,
\[
    \mathcal R_j u_k
    =
    R_jR_{j-1}\cdots R_1u_k
    =
    a_{j,k}u_k .
\]
Hence
\[
    \mathcal R_j f
    =
    \sum_{k=1}^N a_{j,k}\theta_k u_k .
\]
By Theorem~\ref{thm:bias-equals-residual},
\[
    \ibias_j=\|\mathcal R_j f\|^2
    =
    \sum_{k=1}^N a_{j,k}^2\theta_k^2 .
\]
Since \(M_j=I-\mathcal R_j\), we also have
\[
    M_j u_k=(1-a_{j,k})u_k .
\]
Using the variance formula from Theorem~\ref{thm:bias-equals-residual},
\[
    \ivar_j
    =
    \sigma^2\operatorname{tr}(M_jM_j^*)
    =
    \sigma^2
    \sum_{k=1}^N (1-a_{j,k})^2 .
\]

Finally, because \(0\leq \lambda_{i,k}\leq 1\),
\[
    0\leq a_{j+1,k}\leq a_{j,k}\leq 1 .
\]
Thus \(a_{j,k}^2\) is nonincreasing in \(j\), while \((1-a_{j,k})^2\) is nondecreasing in \(j\). Therefore \(\ibias_j\) is nonincreasing and \(\ivar_j\) is nondecreasing. Finally, for nonnegative \(s\) and \(t\) with \(s+t>0\), the quantity \(s/(s+t)\) is nondecreasing in \(s\) and nonincreasing in \(t\); applying this with \(s=\ibias_j\) and \(t=\ivar_j\) proves the monotonicity of~\(\ibr_j\).
\end{proof}

The spectral model also makes the full bias--variance interplay quantitative. The next result gives an exact expression, within this model, for the change in the integrated error produced by adding one multiscale level, and thereby identifies when additional levels stop helping.

\begin{corollary}[Per-level error increment and crossover]\label{cor:mise-increment}
Assume the setting of Theorem~\ref{thm:spectral-damping}. For a fixed mode~$k$, write~$\lambda_{j+1,k}\in[0,1]$ for the spectral response of level~$j+1$, so that the residual factor of that level is~$1-\lambda_{j+1,k}$ and~$a_{j+1,k}=(1- \lambda_{j+1,k})a_{j,k}$. Then the level-$(j+1)$ increment of the integrated error decomposes over modes as
\begin{equation*}
    \begin{split}
    \mise_{j+1}&-\mise_{j}=\sum_{k=1}^N \Delta_{j,k},
    \\
    \Delta_{j,k} &= a_{j,k}\,\lambda_{j+1,k} \cdot \Big[\,2\sigma^2-a_{j,k}(2-\lambda_{j+1,k})(\theta_k^2+\sigma^2)\,\Big].
    \end{split}
\end{equation*}
Consequently, adding level~$j+1$ strictly decreases the error contributed by mode~$k$ (that is,~$\Delta_{j,k}<0$) if and only if~$a_{j,k}\lambda_{j+1,k}>0$ and 
\begin{equation}\label{eq:crossover}
    a_{j,k}\,(2 - \lambda_{j+1,k})>\frac{2\sigma^2}{\theta_k^2+\sigma^2}.
\end{equation}
In particular:
\begin{enumerate}[label=(\roman*)]
\item If~$\sigma=0$, then~$\Delta_{j,k}\le 0$ for all~$j,k$, so~$\mise_j$ is nonincreasing in~$j$: additional levels never increase the error.
\item Since~$a_{j,k}$ is nonincreasing in~$j$ and~$2-\lambda_{j+1,k}\le 2$, once $a_{j,k}<\sigma^2/(\theta_k^2+\sigma^2)$ the criterion~\eqref{eq:crossover} fails at that level and at every subsequent level; that mode's error contribution is then nondecreasing under subsequent levels.
\end{enumerate}
\end{corollary}

\begin{proof}
By Theorem~\ref{thm:spectral-damping}, the per-mode contribution to the integrated error after~$j$ levels is
\[
    m_{j,k}:=a_{j,k}^2\theta_k^2+\sigma^2(1-a_{j,k})^2,\qquad \mise_j=\sum_{k=1}^N m_{j,k}.
\]
Write~$a:=a_{j,k}$ and~$a'=a_{j+1,k}=(1-\lambda_{j+1,k}) a$. Then $(a')^2-a^2=-a^2\cdot\lambda_{j+1,k}\cdot(2-\lambda_{j+1,k})$ and
\[
    (1-a')^2-(1-a)^2=(a-a')\big(2-a-a'\big)=a\cdot\lambda_{j+1,k}\cdot\big(2-a(2-\lambda_{j+1,k})\big).
\]
Hence
\[
    \Delta_{j,k}=\theta_k^2(a'^2-a^2)+\sigma^2\big[(1-a')^2-(1-a)^2\big]
    =a\cdot \lambda_{j+1,k} \cdot\big[2\sigma^2-a(2- \lambda_{j+1,k})(\theta_k^2+\sigma^2)\big],
\]
and summing over~$k$ gives the increment of~$\mise$. If $a_{j,k}\lambda_{j+1,k}=0$, then $\Delta_{j,k}=0$.
Otherwise, its sign equals that of the bracket, which yields
the crossover criterion~\eqref{eq:crossover}. 

For~(i), setting~$\sigma=0$ gives~$\Delta_{j,k}=-a_{j,k}^2\cdot\lambda_{j+1,k}\cdot (2-\lambda_{j+1,k})\theta_k^2\le0$. For (ii), suppose $a_{j,k}<\sigma^2/(\theta_k^2+\sigma^2)$. For every $j'\geq j$, monotonicity of $a_{j',k}$ gives
\[
a_{j',k}(2-\lambda_{j'+1,k})
\leq 2a_{j',k}
\leq 2a_{j,k}
<\frac{2\sigma^2}{\theta_k^2+\sigma^2}.
\]
Thus the criterion remains violated at every subsequent level.
\end{proof}

\begin{remark}[Number of levels: what the spectral model does and does not say]\label{rem:optimal-levels}
Corollary~\ref{cor:mise-increment} exhibits the bias--variance tradeoff at the level of the error: each additional level lowers the bias contribution of every mode but raises its variance contribution, and the net effect on mode~$k$ is governed by~\eqref{eq:crossover}. A mode benefits from an additional level precisely when the strict crossover criterion~\eqref{eq:crossover} holds. Two qualifications are in order before reading this as a stopping rule. First, the criterion concerns the MISE and not the bias ratio: within the model, $\ibr_j$ is nonincreasing for every~$j$ (Theorem~\ref{thm:spectral-damping}), so the bias ratio by itself never signals a stopping level, and the choice of the number of levels is an MSE consideration. Second, the corollary does not single out any particular number of levels. To see how the optimum depends on the parameters, consider a single dominant mode~$k$ with coefficient~$\theta:=\theta_k$ and a stationary filter, $\lambda_{i,k}=\lambda\in(0,1)$ for all~$i$, so that~$a_{j,k}=(1-\lambda)^j$ and
\[
    \mise_j=\theta^2(1-\lambda)^{2j}+\sigma^2\big(1-(1-\lambda)^j\big)^2 .
\]
Viewed as a function of~$\beta=(1-\lambda)^j\in(0,1]$, this is minimized at $\beta^\ast=\sigma^2/(\theta^2+\sigma^2)$; since~$\beta$ decreases with~$j$, the error decreases until~$(1-\lambda)^j$ crosses~$\beta^\ast$ and increases afterward. With $\sigma<0$, the real-valued optimum is
\[
    j^\ast=\frac{\log\!\big(\sigma^2/(\theta^2+\sigma^2)\big)}{\log(1-\lambda)}
    =\frac{\log\!\big(1+\theta^2/\sigma^2\big)}{-\log(1-\lambda)} .
\]
The optimum grows only logarithmically with the signal-to-noise ratio~$\theta^2/\sigma^2$ of the mode, but it grows like~$1/\lambda$ as the spectral response~$\lambda$ decreases, since~$-\log(1-\lambda)\approx\lambda$ for small~$\lambda$. For instance, for~$\theta^2/\sigma^2=100$ one obtains~$j^\ast\approx2.0$ for~$\lambda=0.9$, $j^\ast\approx6.7$ for~$\lambda=0.5$, and~$j^\ast\approx44$ for~$\lambda=0.1$. Thus, the optimal number of levels is small only when the level operators damp the relevant modes aggressively; when they act weakly on a mode, many levels may be needed before the crossover is reached, even at moderate signal-to-noise ratios. In the general multimode case, different modes cross over at different levels, and the optimal number of levels is determined jointly by the energy distribution~$(\theta_k^2)_k$ and by the spectral responses~$(\lambda_{i,k})$; modes that are reproduced exactly by every level, such as the constant mode when~$\lambda_{i,1}=\phi_i(0)=1$, have~$a_{j,1}=0$ for~$j\geq1$ and are unaffected by further levels. A second, noise-independent mechanism limits the marginal benefit of additional levels for the deterministic residual itself: in the multilevel MLS analysis of~\cite{durst2026improved}, the error keeps decreasing with the number of levels, but the acceleration of the convergence order is not uniform in the level and consists of an initial boost over the first few levels followed by a plateau \cite[Lemma~3.5 and the discussion following it]{durst2026improved}; in their experiments, the observed rates stagnate after roughly six levels. The three-level choice in Section~\ref{sec:num-settings} is therefore an empirical choice for the operators, hierarchies, and noise levels considered there, informed by these mechanisms but not implied by them. In practice, the number of levels can be selected by monitoring the estimated MSE across levels or, when the noise level is known, through the criterion~\eqref{eq:crossover} applied to the dominant modes.
\end{remark}

The preceding theorem does not require the exact assumptions to hold in every numerical example. Rather, it identifies the mechanism: multiscale error correction replaces the original bias by a product of residual filters. The experiments below test how robustly this mechanism appears for MLS, Shepard approximation, and manifold-valued data.

\subsection{Residual-bias interpretation in manifold settings}\label{sec:manifold}

We give a local, conditional extension of the scalar theory. Fix the data sites and the hierarchy, write $q=\dim\M$, and fix a normal coordinate chart centered at a point $p\in\M$. Throughout this section, $C$ denotes a generic constant that may depend on the fixed number of levels, on the data size, and on uniform stability and geometric bounds, but not on $r$ or on $\sigma$. The three assumptions below are collected explicitly because they are what the lemmas of this section use; they are sufficient rather than necessary.

\begin{assumption}[Local geometry]\label{ass:geometry}
The target values, the noisy data, all intermediate manifold values, and the outputs of the scheme lie in a strongly convex normal neighborhood of $p$ of radius $Cr$, on which the sectional curvature and the first two derivatives of the chart are uniformly bounded. Moreover, $r,\sigma\leq1$ are sufficiently small.
\end{assumption}

Write $u(x)=\log_pF(x)$ and $u_X=(u(x_i))_{i=1}^N\in(T_p\M)^N$. Let $A_n(x):(T_p\M)^N\to T_p\M$ be the deterministic linear map obtained by performing the corresponding Euclidean multiscale scheme and evaluating at $x$. 

Assume that $\sup_{x\in\Omega}\|A_n(x)\|_{\mathrm{op}}\leq C$,
uniformly in the local parameters $r$ and $\sigma$ and define the deterministic evaluation residual
\[
    e_n(x)=u(x)-A_n(x)u_X.
\]
On the data sites, this is the vector-valued residual $((I-M_n)\otimes I_q)u_X$; on an external evaluation set it is an approximation residual, not a square data-site operator applied to $u$.

\begin{assumption}[Local consistency of the recursion]\label{ass:consistency}
The full manifold recursion satisfies, almost surely and uniformly in $x$,
\begin{equation}\label{eq:full-local-consistency}
 w(x):=\log_p\widetilde F_n(x)
      =A_n(x)u_X^\varepsilon+g_n(x),
 \quad \|g_n(x)\|\leq Cr^3 .
\end{equation}
\end{assumption}

Note that the consistency of an individual intrinsic average does not by itself imply \eqref{eq:full-local-consistency}: residual formation, tangent-vector identification, and the update must also be consistent. A sufficient condition is that each primitive of the fixed-length recursion agrees with its Euclidean counterpart up to $O(r^3)$ in these coordinates, with uniformly bounded local Lipschitz constants (for example, for the action of elementary arithmetic operators). Indeed, propagating the errors through the finite composition gives $c_{j+1}\leq C_jc_j+C_jr^3$, and hence $c_j=O(r^3)$. For intrinsic means, the required expansion is
\[
 \log_p\operatorname{Av}_{\M}(\{\exp_pv_i\},\{a_i\})
 =\sum_i a_iv_i+O(r^3),
 \qquad \sum_i a_i=1,
\]
under the usual local existence and stability assumptions. For nonnegative weights summing to one, this expansion indeed follows from the local barycenter expansion, as seen in \cite[Theorem~5]{pennec2019curvature}, applied to the discrete probability measure $\sum_i a_i\delta_{\exp_p(v_i)}$ (the model used there to represent the points and weights). The third order stems from the curvature correction. The consistency estimate includes any errors introduced
by inexact averaging steps.

\begin{assumption}[Coordinate noise]\label{ass:noise}
The noise in normal coordinates has bounded local support,
\begin{equation}\label{eq:bounded-noise}
 \|u_X^\varepsilon-u_X\|\leq Cr \qquad\text{almost surely},
\end{equation}
and satisfies
\begin{equation}\label{eq:coordinate-noise-local}
 \mathbb E u_X^\varepsilon=u_X+\beta,
 \quad \|\beta\|\leq C\sigma^2,
 \quad \zeta=u_X^\varepsilon-\mathbb E u_X^\varepsilon,
 \quad \Sigma=\mathbb E[\zeta\zeta^*],
 \quad \operatorname{tr}\Sigma\leq C\sigma^2.
\end{equation}
\end{assumption}

Assumption~\ref{ass:noise} allows anisotropic and correlated noise. For example, centered, bounded small tangent perturbations $F_i^\varepsilon=\exp_{F(x_i)}\xi_i$ with $\|\xi_i\|\leq Cr$ almost surely and $\mathbb E\|\xi_i\|^2=O(\sigma^2)$ satisfy these estimates when the coordinate maps and their first two derivatives are uniformly bounded on the perturbation paths. Taylor expansion of $v\mapsto\log_p\exp_{F(x_i)}v$ proves the mean estimate because its linear term has zero expectation. Its local Lipschitz bound proves the second-moment estimate.

The locality and noise assumptions are imposed jointly: bounded input perturbations alone do not ensure that every intermediate iterate remains in the prescribed neighborhood. The bounded output support in Assumption~\ref{ass:geometry} justifies differentiation under the expectation below. The untruncated Gaussian experiments illustrate the local mechanism but do not satisfy the bounded-support hypotheses.

\begin{lemma}[Local tangent-space bias]\label{lem:manifold-linearization}
Under Assumptions~\ref{ass:geometry}--\ref{ass:noise}, assume the intrinsic mean $\overline F_n(x)=\mathbb E_{\M}\widetilde F_n(x)$ exists uniquely in the same neighborhood. For fixed $n$, with $\chi=r^3+\sigma^2$,
\[
 \log_p\overline F_n(x)=A_n(x)u_X+O(\chi),
\]
and
\[
 \bias[\widetilde F_n,F](x)=\|e_n(x)\|^2+
 O\bigl(r^2\|e_n(x)\|^2+\chi\|e_n(x)\|+\chi^2\bigr).
\]
\end{lemma}

\begin{proof}
Put $b=\log_p\overline F_n(x)$ and suppress $x$ temporarily. The normal-coordinate expansion of the gradient of squared distance, uniformly for $\|b\|,\|w\|=O(r)$, gives
\[
 \nabla_b\tfrac12\rho(\exp_pb,\exp_pw)^2=b-w+O(r^3).
\]
The first-order condition for the intrinsic mean, with differentiation under the expectation justified by the bounded output support in
Assumption~\ref{ass:geometry}, yields $b=\mathbb Ew+O(r^3)$. By \eqref{eq:full-local-consistency} and \eqref{eq:coordinate-noise-local},
\[
 \mathbb Ew=A_n(x)(u_X+\beta)+\mathbb Eg_n(x)
           =A_n(x)u_X+O(\chi).
\]
Uniform stability of $A_n(x)$ is used here. Finally, the metric in normal coordinates satisfies
\[
 \rho(\exp_pa,\exp_pb)^2
 =\|a-b\|^2+O(r^2\|a-b\|^2).
\]
Apply this with $a=u(x)$ and $b=A_n(x)u_X+O(\chi)$ and expand the squared norm. Since $r\leq1$, the smaller mixed remainders are absorbed into the displayed bound.
\end{proof}

Let $\mu$ be a finite measure on the evaluation domain $\Omega$, and assume that all preceding estimates hold uniformly there in the fixed chart. Set
\[
 \mathcal B_n^2=\int_\Omega\|e_n(x)\|^2\,d\mu(x),
 \qquad
 \mathcal V_n=\int_\Omega
       \operatorname{tr}\bigl(A_n(x)\Sigma A_n(x)^*\bigr)\,d\mu(x).
\]
These definitions also cover a finite evaluation grid with quadrature weights. They do not identify an external-grid trace with a data-site trace.

\begin{lemma}[Integrated tangent-space variance and mean squared error]
\label{lem:manifold-variance}
Under the hypotheses of Lemma~\ref{lem:manifold-linearization}, with
$\chi=r^3+\sigma^2$,
\[
 \ivar_n=\mathcal V_n+O(r^2\sigma^2+r^3\sigma+r^6),
\]
\[
 \mise_n=\mathcal B_n^2+\mathcal V_n+
 O\bigl(r^2\mathcal B_n^2+\chi \mathcal B_n+r^2\sigma^2+\chi\sigma+\chi^2\bigr).
\]
If, additionally, $\Sigma=\sigma^2 I_{Nq}+O(r\sigma^2)$ in operator norm, then
\[
 \mathcal V_n=\sigma^2\int_\Omega\operatorname{tr}
       \bigl(A_n(x)A_n(x)^*\bigr)\,d\mu(x)+O(r\sigma^2).
\]
In particular, if $A_n(x)=a_n(x)\otimes I_q$ for a scalar row operator $a_n(x)$, the leading isotropic term is $q\sigma^2\int_\Omega\|a_n(x)\|_2^2\,d\mu(x)$, where $\sigma^2$ is the variance per tangent coordinate.
\end{lemma}

\begin{proof}
Fix $x$, write $A=A_n(x)$, $g=g_n(x)$, and use the notation of the previous proof. Since $b-\mathbb Ew=O(r^3)$,
\[
 w-b=A\zeta+\vartheta,\qquad
 \vartheta=g-\mathbb Eg+(\mathbb Ew-b),\qquad \|\vartheta\|\leq Cr^3.
\]
Consequently, by Cauchy--Schwarz and $\mathbb E\|A\zeta\|^2=\operatorname{tr}(A\Sigma A^*)=O(\sigma^2)$,
\[
 \mathbb E\|w-b\|^2
 =\operatorname{tr}(A\Sigma A^*)+O(r^3\sigma+r^6).
\]
The normal-coordinate distance estimate adds $O(r^2\mathbb E\|w-b\|^2)$, giving the pointwise variance formula.

For the error, write $k=A\beta+g$, so $\|k\|\leq C\chi$ and
\[
 u(x)-w=e_n(x)-A\zeta-k.
\]
The cross term between the deterministic residual and $A\zeta$ has zero expectation. The other cross terms are bounded by $C\chi\|e_n(x)\|$ and $C\chi\sigma$, respectively. Thus
\[
 \mathbb E\|u(x)-w\|^2
 =\|e_n(x)\|^2+\operatorname{tr}(A\Sigma A^*)
 +O\bigl(\chi\|e_n(x)\|+\chi\sigma+\chi^2\bigr).
\]
The metric correction contributes $O(r^2\|e_n(x)\|^2+r^2\sigma^2+r^2\chi^2)$, with the last term absorbed into $O(\chi^2)$. Integrate these estimates and use
\[
 \int_\Omega\|e_n(x)\|\,d\mu(x)
 \leq\mu(\Omega)^{1/2}\mathcal B_n.
\]
For isotropy, write $\Sigma=\sigma^2 I+H$ and use $|\operatorname{tr}(AHA^*)|\leq\|H\|_{\mathrm{op}}\|A\|_F^2$. Finally $\|a_n(x)\otimes I_q\|_F^2=q\|a_n(x)\|_2^2$.
\end{proof}

In Euclidean space, for an exact linear implementation and centered additive noise, $g_n=0$ and $\beta=0$, and the distance corrections vanish. The exact vector-valued identities are then $\ivar_n=\mathcal V_n$ and $\mise_n=\mathcal B_n^2+\mathcal V_n$. On a curved manifold the preceding lemmas give absolute error estimates. Their interpretation as a relative approximation to the bias ratio requires the remainders to be small compared with the MSE; it is not uniform when the denominator tends to zero.

We finally describe the scalar Shepard matrices underlying this local linear model. Let $X_i\subseteq X$ consist of distinct sites, let $\Pi_i:\mathbb R^X\to\mathbb R^{X_i}$ be restriction, and define
\[
 (B_i)_{aj}=
 \frac{K(\|x_a-x_j^{(i)}\|/\delta_i)}
 {\sum_{\ell}K(\|x_a-x_\ell^{(i)}\|/\delta_i)},
 \qquad Q_i=B_i\Pi_i.
\]
Here $K$ is the radial kernel of the scheme and $\delta_i$ its level-$i$ support radius, as in Section~\ref{sec:mls}. Assume the kernel values are nonnegative and every denominator is strictly positive. Write
\[
 (\widehat K_i)_{j\ell}
 =K(\|x_j^{(i)}-x_\ell^{(i)}\|/\delta_i),
 \qquad \widehat D_i=\operatorname{diag}(\widehat K_i\mathbf1).
\]

\begin{lemma}[Structure of the Shepard level operators]
\label{lem:shepard-structure}
Assume $\widehat K_i$ is symmetric positive definite, in addition to the preceding positivity and coverage assumptions. Then $Q_i$ is nonnegative and reproduces constants. Its nonzero eigenvalues, including algebraic multiplicities, are those of $\Pi_iB_i=\widehat D_i^{-1}\widehat K_i$ and belong to $(0,1]$. The eigenvalue zero has algebraic and geometric multiplicity $|X|-|X_i|$. In particular, $\operatorname{spec}(Q_i)\subset[0,1]$ and the spectral radius of~$Q_i$ equals one.
\end{lemma}

\begin{proof}
The matrices $B_i$ and $Q_i$ are row-stochastic, so $Q_i\mathbf1=\mathbf1$ and $\|Q_i\|_\infty=1$. Set $N=|X|$, $N_i=|X_i|$, and $S_i=\Pi_iB_i$. The matrix $S_i$ is similar to
\[
 \widehat D_i^{1/2}S_i\widehat D_i^{-1/2}
 =\widehat D_i^{-1/2}\widehat K_i\widehat D_i^{-1/2},
\]
which is symmetric positive definite. Hence its eigenvalues are real and positive; row-stochasticity bounds them above by one. The identity
\[
 \det(tI_N-B_i\Pi_i)=t^{N-N_i}\det(tI_{N_i}-\Pi_iB_i)
\]
proves the eigenvalue assertions and the algebraic multiplicity of zero. Since $S_i$ is invertible, $B_i$ has full column rank and $\ker Q_i=\ker\Pi_i$, of dimension $N-N_i$. Thus the geometric multiplicity is the same.
\end{proof}

The spectrum of each individual level does not imply Euclidean norm contraction of the residual product. The $Q_i$ need not be normal in a common inner product or commute across levels. Thus the joint spectral theorem does not apply to general subset-based Shepard schemes. The local bias identities remain available under the full-recursion consistency assumption, but actual bias reduction additionally requires control of the deterministic residual. Such control is a hypothesis or an empirical observation for the present schemes, not a consequence of the location of their individual spectra. For MLS, the exact scalar linearity argument still applies, whereas the nonnegative-weight Shepard lemma generally does not.

\section{Numerical evaluation}\label{sec:numerics}
We numerically demonstrate the effect of the multiscale scheme on the bias ratio~\eqref{eq:bias-ratio}. The results in Section~\ref{sec:residual-bias} identify residual contraction as a sufficient mechanism for bias reduction and provide a local tangent-space explanation for manifold-valued data. For the scalar MLS experiments, the residual-bias identity of Theorem~\ref{thm:bias-equals-residual} holds exactly for every realized design; since the data sites and the hierarchy are resampled in every trial, the reported bias is the deterministic multiscale residual averaged over designs, and the reported variance carries an additional design component, as quantified in Proposition~\ref{prop:random-design}. The per-level contraction of Corollary~\ref{cor:residual-contraction} is the hypothesis under test (Remark~\ref{rem:contraction-hypothesis}); for the manifold Shepard experiments, the same first-tier route applies to leading order through Lemmas~\ref{lem:manifold-linearization}--\ref{lem:manifold-variance}, while Lemma~\ref{lem:shepard-structure} shows that the per-level operators are constant-reproducing averaging operators with real spectrum in~$[0,1]$ but do not satisfy the joint hypotheses of the spectral model. The following demonstrations test how robustly the mechanism appears in practical approximation methods, target functions, and noise levels. The final experiment, with diffusion operators, is designed so that the hypotheses of the spectral model hold exactly, and it compares the measured integrated quantities with the theoretical predictions.

We begin the numerical experiments by applying the multiscale method to moving least squares (MLS) approximation with different polynomial degrees. The target function in this setting is a smooth scalar-valued function, and we perform experiments across datasets of varying sizes. The results show that applying the multiscale method achieves performance comparable to increasing the MLS polynomial degree by one, while also achieving a substantially lower bias ratio and improved error rates. Next, we turn to Shepard approximation~\cite{shepard1968two} and extend the setting from Euclidean to manifolds. In particular, we compare the standard Shepard approximation with its multiscale counterpart for a noisy, high-bias SO(3)-valued target function, analyzing the performance as a function of the signal-to-noise ratio (SNR). In addition, we extend the experiments to another manifold-valued setting by considering target functions on the manifold of symmetric positive definite (SPD) matrices. Finally, using a diffusion operator, we demonstrate the theoretical properties of the bias ratio established in Section~\ref{sec:spectral}.

\subsection{Settings of the numerical framework}\label{sec:num-settings}

The following settings of the numerical experiments apply to the MLS and manifold-valued experiments; the fixed-grid diffusion experiment is described separately in Section~\ref{sec:diffusion}.

We use a scattered dataset~$X$, which is sampled uniformly from~$[0,1]^2$. Three parameters determine the hierarchical sequence of nested subsets: 1) the dataset~$X$, 2) the number of levels, and 3) the proportional growth rate~$\gamma$. The growth rate~$\gamma$ controls the size of each subset in the hierarchical sequence; we set~$\gamma=80\%$. All the results are for a three-level multiscale. In preliminary runs with the operators, noise levels, and hierarchies considered here, the marginal contribution of further levels was small; we stress that this is an empirical choice for the present setting rather than a general rule, and we refer to Remark~\ref{rem:optimal-levels} for the mechanisms that govern the beneficial number of levels. We generate the hierarchy by first sampling $X_n=X$, then sampling $X_{n-1}$ as a random subset of $X_n$, etc., ensuring $ X_1 \subset \ldots \subset X_n = X$. For each data point in the graphs, we repeat the process of generating datasets, hierarchical sequences, and sample noise, and apply approximations for~$100$ trials. To evaluate the experimental error, we generate a dense, uniformly spaced grid over $[0,1]^2$, where we compare the target function and its different approximations.

The design is therefore resampled in every trial, independently of the noise: the data sites, the nested subsets, and the support radii $\delta_i=\nu h_i$, which are recomputed from the realized mesh norm of each level. The reported quantities average over both the design and the
observational noise. For scalar linear estimators, the resulting
design contribution is described by
Proposition~\ref{prop:random-design}; under heteroscedastic noise,
its noise-variance term is replaced by
$\mathbb E_{\mathbf X}
[a_n(x;\mathbf X)\Sigma_{\mathbf X}a_n(x;\mathbf X)^*]$,
where $\Sigma_{\mathbf X}
=\operatorname{Cov}(\varepsilon\mid\mathbf X)$.
For manifold-valued estimators, the reported quantities use
the intrinsic definitions of Section~\ref{sec:bvt}. The diffusion experiment of Section~\ref{sec:diffusion} is the only one in which the design and the level operators are held fixed across trials, as assumed in Theorem~\ref{thm:spectral-damping}.

In the numerical results, we track two key metrics: the MSE~\eqref{eq:mse} (and its manifold analogue~\eqref{eq:mse-manifold}) and the bias ratio~\eqref{eq:bias-ratio}. Both are pointwise quantities, estimated at each point of the evaluation grid from the $100$ trials, with the average over trials in~\eqref{eq:mc-mean} (respectively, the intrinsic mean in the manifold case) playing the role of $\E[\widetilde f(x)]$. Each curve then reports the distribution of the resulting pointwise quantity over the evaluation grid: the solid line is its median over the evaluation points, and the shaded band is the 25th--75th percentile range over the evaluation points. The bands therefore measure the spatial variability of the pointwise error over the domain, and not Monte Carlo variability across trials, which has already been averaged out in forming the pointwise bias and variance. In each experiment, we compare a chosen approximation scheme (referred to as ``singlescale'') and its multiscale counterpart. The code is publicly available\footnote{\url{https://github.com/ABASASA/Multiscale_Bias_Ratio}}.

\subsection{Multiscale based on moving least squares (MLS) over scattered data} \label{sec:mls}

We examine the multiscale for moving least squares approximation (MLS). Our MLS implementation is based on inverting the system of equations, and we compare the behavior of MLS with polynomial degrees of 1, 2, and 3. We select a Gaussian radial kernel, $K(r)=e^{-2r^2}$, evaluated at the normalized distance $r=\delta_i^{-1}\|x-s\|_2$ between an evaluation point~$s$ and a data site~$x$; the same convention is used for the Shepard kernel in Section~\ref{sec:so3}. The support parameter is set per level as~$\delta_i=\nu h_i$ with~$\nu=3$ and~$h_i$ the mesh norm (fill distance) of the level-$i$ data set~$X_i$; in the terminology of~\cite{durst2026improved}, $\nu$ is the ratio between the support radius and the mesh width.

As an initial experiment and preliminary step toward manifold-valued functions, we compare the MLS and its multiscale behaviors over a smooth, scalar-valued target function,
\begin{equation}\label{eq:mls-target}
f(x,y) = e^{x^2 +  y^2} + 3 .
\end{equation}
Without added noise, this target function, being analytic, admits a high-order approximation. In this noiseless experiment, all randomness comes from resampling the data sites and the hierarchy, rather than from observational noise. By Proposition~\ref{prop:random-design}(iv) with $\sigma=0$, the whole MSE is then the mean squared deterministic residual $\E_{\mathbf X}[e_n(x;\mathbf X)^2]$, and the bias ratio measures the fraction of that residual which is systematic across designs; in particular it is strictly below one whenever the residual depends on the realized design, and a lower bias ratio means that a larger share of the remaining error is design dependent rather than reproducible. We investigate the multiscale behavior over various dataset sizes and MLS's approximation degree.

\begin{figure}[!ht]
\centering
\begin{subfigure}{.32\textwidth}
\centering
\includegraphics[width=\textwidth]{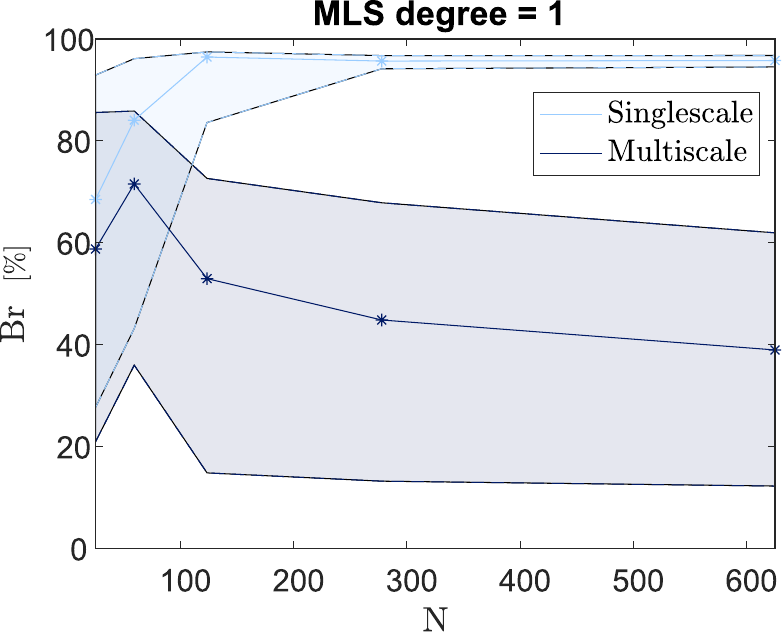}
\caption{Degree 1}\label{fig:mls-br-d1}
\end{subfigure}
\hfill
\begin{subfigure}{.32\textwidth}
\centering
\includegraphics[width=\textwidth]{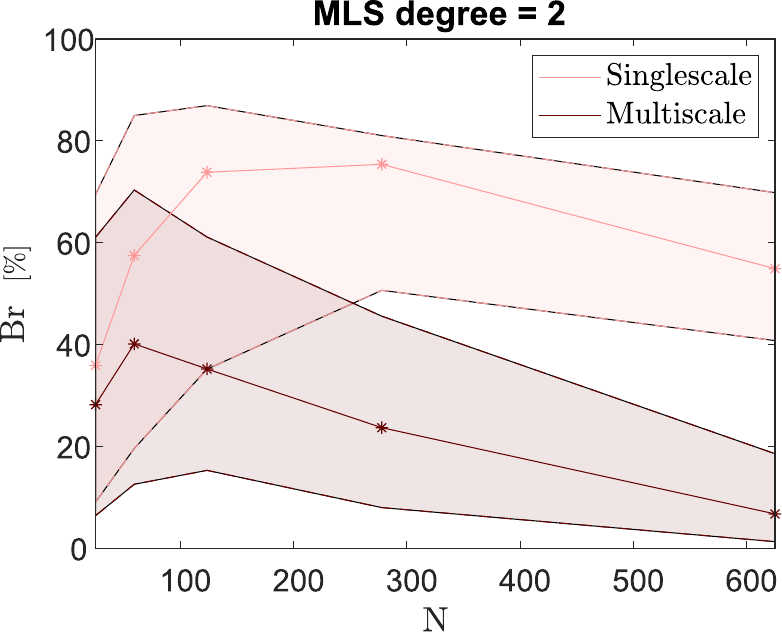}
\caption{Degree 2}\label{fig:mls-br-d2}
\end{subfigure}
\hfill
\begin{subfigure}{.32\textwidth}
\centering
\includegraphics[width=\textwidth]{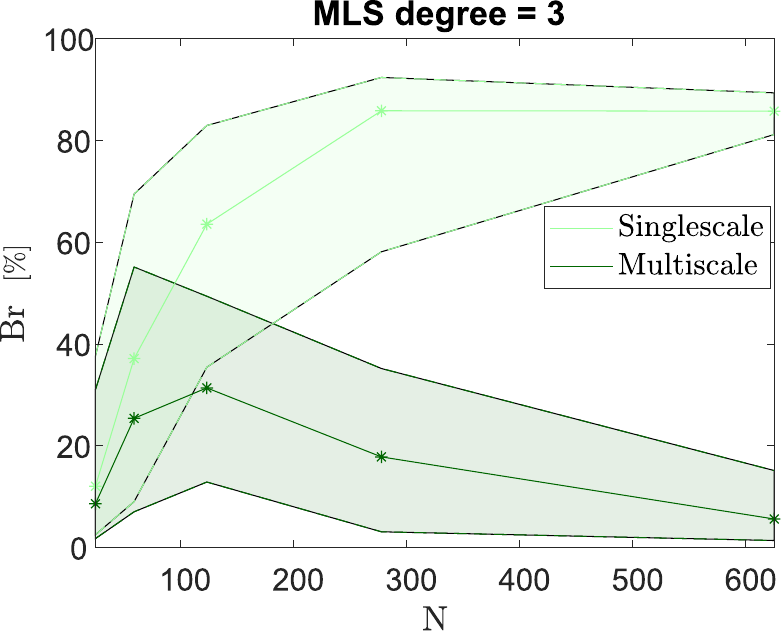}
\caption{Degree 3}\label{fig:mls-br-d3}
\end{subfigure}
\caption{Bias ratio for MLS and multiscale MLS across dataset sizes ($N$) and polynomial degrees. The MLS graphs are shown in lighter colors, and their multiscale counterparts in darker colors: linear (blue), quadratic (red), and cubic (green). The target function is a smooth scalar-valued function~\eqref{eq:mls-target}. Solid line: median over the evaluation grid; shaded region: 25th--75th percentile over the evaluation grid; both are computed from the pointwise quantities estimated over 100 trials. $\br$ is shown in percent.}
\label{fig:mls-br}
\end{figure}

\begin{figure}[!ht]
\centering
\includegraphics[width=0.55\textwidth]{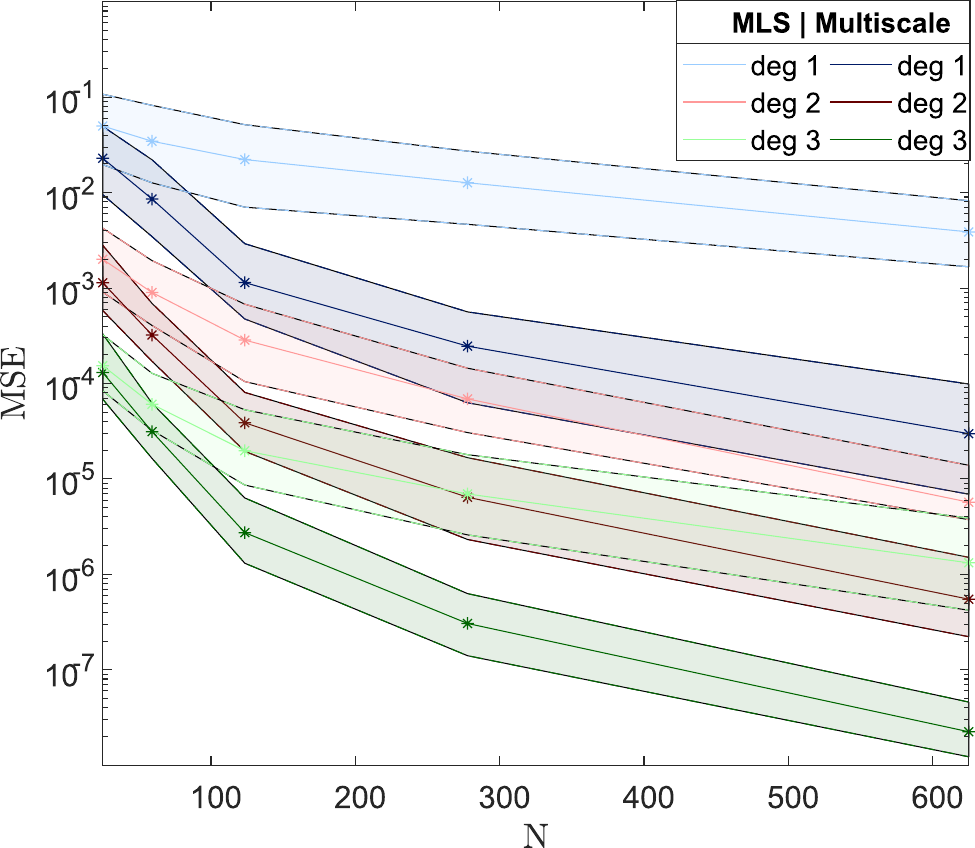}
     \centering   
\caption{MSE for MLS and multiscale MLS across dataset sizes ($N$) and polynomial degrees. The MLS graphs appear in lighter colors, and their multiscale counterparts appear in darker colors: linear (blue), quadratic (red), and cubic (green). The target function is a smooth scalar-valued function~\eqref{eq:mls-target}. In each error curve, the solid line denotes the median over the evaluation grid and the shaded region shows the 25th--75th percentile over the evaluation grid, computed from the pointwise MSE estimated over 100 trials.}      \label{fig:mls-mse}
\end{figure}

Figures~\ref{fig:mls-br} and~\ref{fig:mls-mse} highlight several key observations. In all experiments, the multiscale method improves both the MSE and the bias ratio. In particular, Figure~\ref{fig:mls-br} shows that the multiscale approach consistently reduces the bias ratio across all polynomial degrees and dataset sizes~$N$. Moreover, Figure~\ref{fig:mls-mse} demonstrates that applying the multiscale method is comparable to increasing the degree of the MLS by one. This observation is in line with the convergence analysis of multilevel MLS in~\cite{durst2026improved}: for gridded data and sufficiently smooth targets, error correction over~$L$ levels raises the convergence order of degree-$r$ MLS by up to nearly~$(1-1/L)(r+1)$, that is, by roughly one order for three levels and~$r=1$. Since our sites are random and the weight is Gaussian, the comparison is qualitative. Given the properties of the target function, we expect this behavior to persist for any polynomial degree, except in cases affected by rounding errors, ill-conditioned matrices, or insufficient dataset sizes.

Furthermore, we investigate the multiscale for a noisy, smooth and scalar-valued function. Specifically, we explore the MLS of degree~$1$ for~$f(x,y)$ from~\eqref{eq:mls-target} in noisy settings. We add the noise in a proportional manner, i.e., $f(x,y) \cdot( 1 +\sigma \cdot \epsilon)$, where~$\epsilon$ is a standard normal random variable and~$\sigma$ is given. This noise is heteroscedastic: conditionally on the design, $\operatorname{Cov}(\varepsilon)=\sigma^2\operatorname{diag}\bigl(f(x_i)^2\bigr)$ rather than $\sigma^2I$. The bias identities of Theorem~\ref{thm:bias-equals-residual} and Proposition~\ref{prop:random-design}(ii) require only $\mathbb E[\varepsilon\mid\mathbf X]=0$ and therefore still apply; the variance formulas, which assume $\operatorname{Cov}(\varepsilon)=\sigma^2I$, do not, and the variance is estimated directly. We estimate the SNR level for a given~$\sigma$ via~$\operatorname{SNR} = \nicefrac{1}{\sigma}$. Here, the dataset size is fixed at approximately $21^2$ points, sampled uniformly.

\begin{figure}[!ht]
\centering
\begin{subfigure}{.45\textwidth}
\centering
\includegraphics[width=\textwidth]{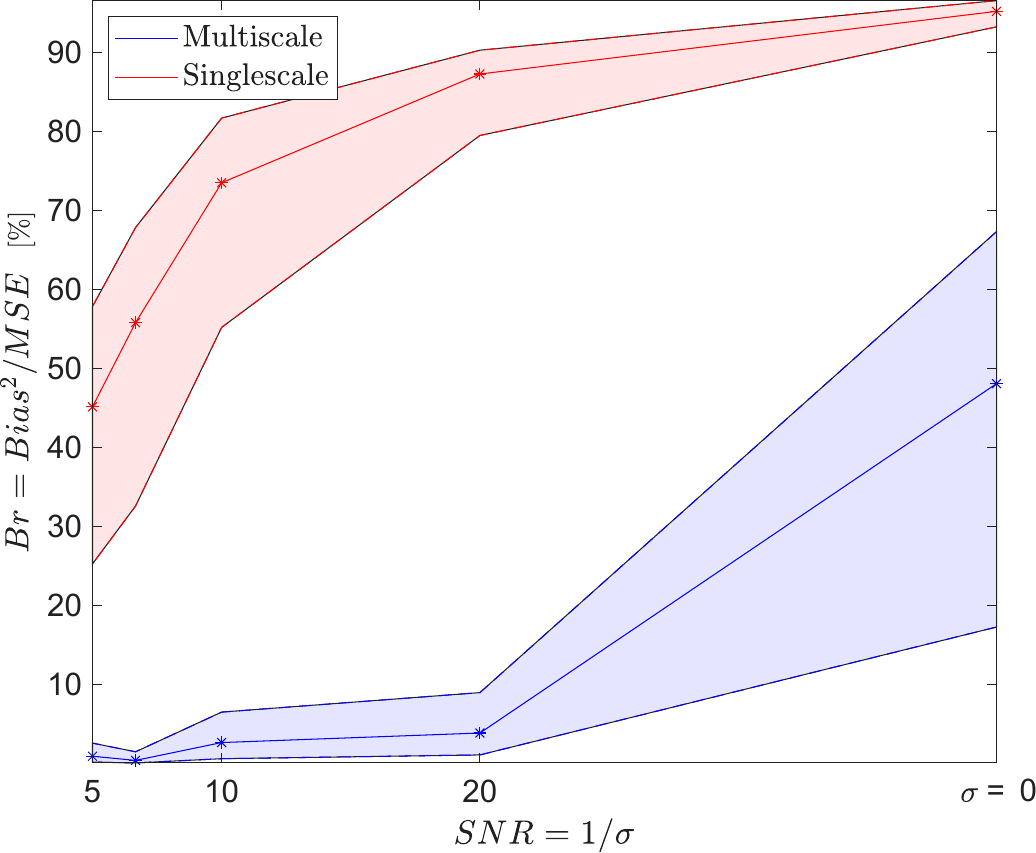} 
\centering
\caption{Bias ratio}\label{fig:mls-noisy-br}
\end{subfigure} \quad
\begin{subfigure}{.45\textwidth}
\includegraphics[width=\textwidth]{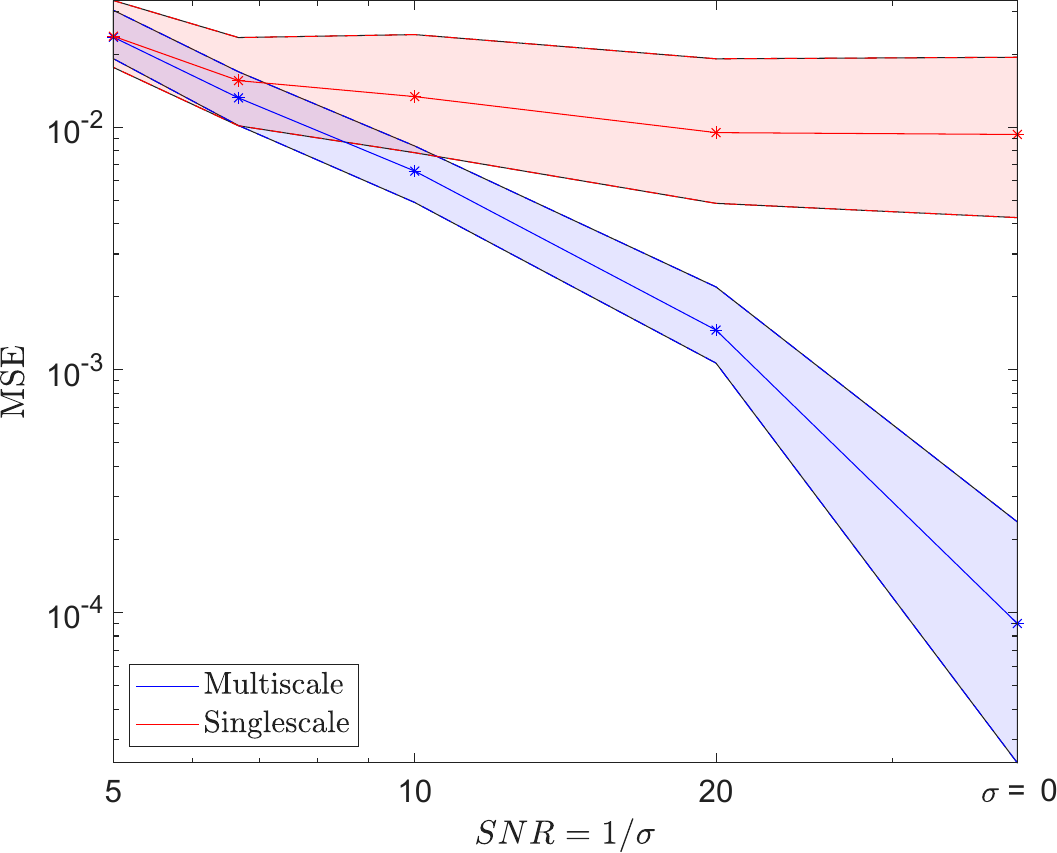}
     \centering   
    \caption{Error rates}
    \label{fig:mls-noisy-mse}
\end{subfigure}%
\caption{Comparing linear MLS (red) and its multiscale (blue) across different SNR (SNR $=\sfrac{1}{ \sigma}$) levels. The target function is the scalar-valued function~\eqref{eq:mls-target}, and the samples are corrupted by multiplicative Gaussian noise, $f(x,y)(1+\sigma\epsilon)$ with $\epsilon$ standard normal. Solid lines denote the median over the evaluation grid and shaded regions show the 25th--75th percentile over the evaluation grid, computed from the pointwise quantities estimated over 100 trials.}
\label{fig:mls-noisy}
\end{figure}

The results, in terms of MSE and bias ratio, are presented in Figure~\ref{fig:mls-noisy}. As can be seen, the multiscale achieves superior performance in terms of MSE and low bias ratio. Note that, in this setting, a higher-degree MLS results in an extremely low bias ratio, $1$--$3\%$, and the MLS performs worse than linear MLS in terms of MSE. In that case, applying multiscale does not improve the results.

\subsection{Approximating a function with values on the SO(3) manifold}\label{sec:so3}

We extend our previous experiments to the rotation manifold setting by introducing noise into the measurements. The Lie group SO(3) is a compact Riemannian manifold, and we follow~\cite{sharon2022multiscale} to apply similar geodesic-based computations, using the Riemannian exponential and the principal logarithm on SO(3).

Consider the smooth target function,
\begin{equation}\label{eq:so3-target}
\begin{split} 
&G:[0,1]^2 \to SO(3), \\
&G(t,w) = \Rot_z\big(1.25\sin(5t)-0.1\big)\, \Rot_y\big(\tfrac{1}{2}w^2 - \sin(3t)\big)\, \Rot_x\big(1.5\cos(2t)\big),
\end{split}
\end{equation}
where~$\Rot_i$ denotes the rotation matrix about the $i$-th axis (e.g., $\Rot_y$ corresponds to rotation about the pitch axis); this notation avoids a clash with the residual operators~$R_i$ of Section~\ref{sec:linear}.

Gaussian noise is introduced in the tangent space at each measurement and then projected back onto the manifold. Performance is assessed using the bias ratio and MSE as functions of the noise standard deviation~$\sigma$. For convenience in plotting, we report results as a function of SNR$\,\coloneqq \nicefrac{1}{\sigma}$. This aligns with the standard SNR definition, since the norm of any rotation matrix is constant, so the ``signal'' part of the SNR is always fixed. The dataset comprises approximately~$13^2$ points, sampled uniformly at random from the domain~$[0,1]^2$.

Following~\cite{sharon2022multiscale, manton2004globally}, we employ an iterative variant of the Shepard approximation, originally introduced in~\cite{shepard1968two}. Note that the manifold experiments illustrate the local theory under Gaussian tangent noise, which does not satisfy the almost sure bound~\eqref{eq:bounded-noise} of Assumption~\ref{ass:noise}. Given a dataset~$X = \{(x_i, F_i) \}_{i=1}^N$, the Shepard approximation at a point~$s$ is computed as follows:
\begin{equation*}
   Y_j = \operatorname{Exp}_{Y_{j-1}}\left( \frac{\sum_{i=1}^N K(r_i) \, \operatorname{Log}_{Y_{j-1}}(F_i)} {\sum_{i=1}^N K(r_i)}\right),
\end{equation*}
where $r_i = \delta_i^{-1} \|x_i - s\|_2$ and the level-dependent support radius~$\delta_i$ is defined as in Section~\ref{sec:mls}. Note that in the product in the numerator, $\operatorname{Log}_{Y_{j-1}}(F_i)$ is a matrix and $K(r_i)$ is a scalar. We use the radial kernel~$K(r) = (\max\{1-r, 0 \})^4 \cdot (4r+1)$ from~\cite{wendland2004scattered}. The iterative process is terminated when either (i) the number of iterations reaches 10, or (ii) the update is within 0.01 in geodesic distance of the previous iterate.

\begin{figure}[!ht]
\centering
\begin{subfigure}{.45\textwidth}
\centering
\includegraphics[width=\textwidth]{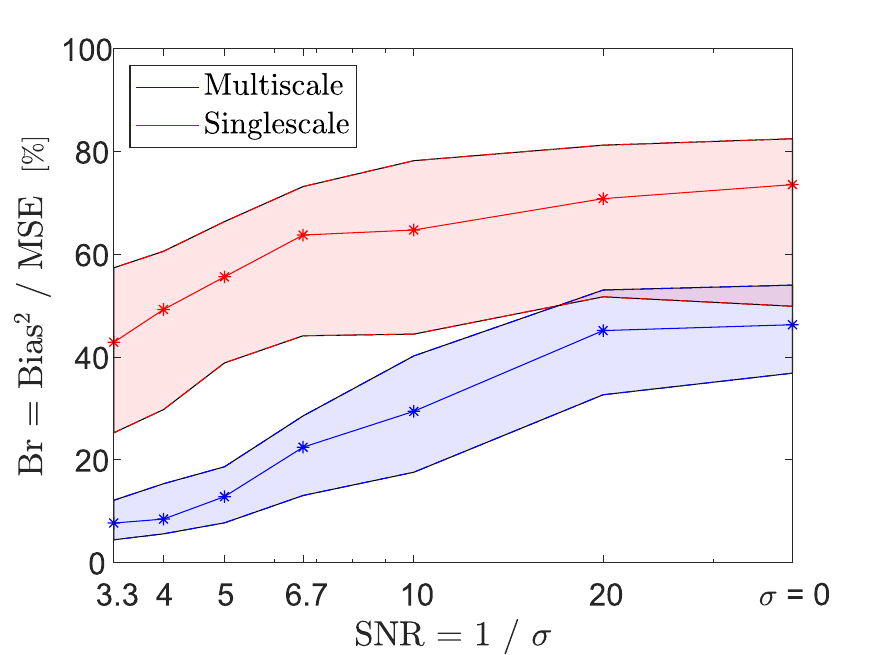} 
\centering
\caption{Bias ratio}\label{fig:so3-br}
\end{subfigure} \quad
\begin{subfigure}{.45\textwidth}
\includegraphics[width=\textwidth]{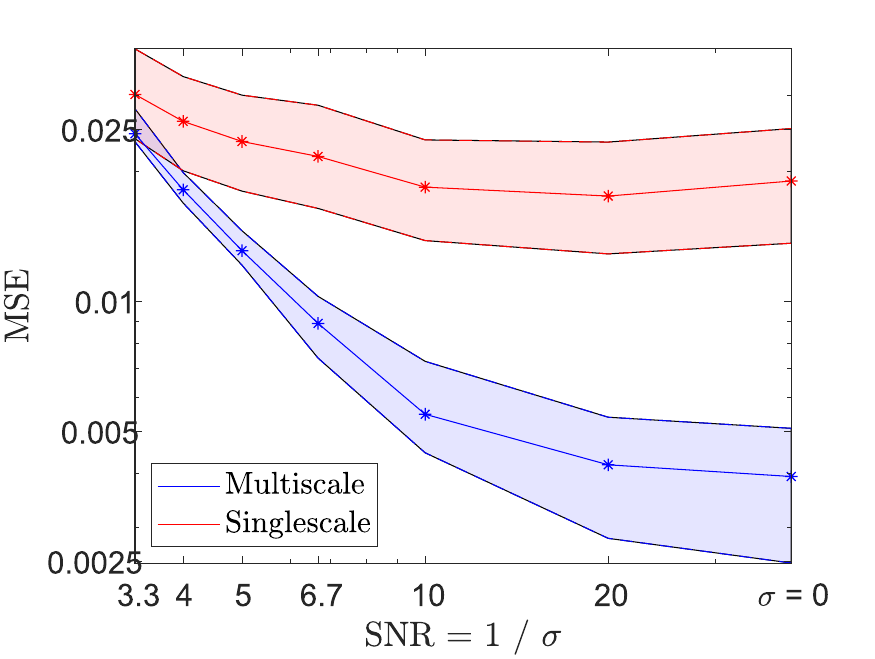}
     \centering   
    \caption{Error rates}
    \label{fig:so3-mse}
\end{subfigure}%
\caption{Comparing Shepard (red) and its multiscale (blue) across different SNR (SNR $=\sfrac{1}{\sigma}$) levels. The target function is the $\operatorname{SO}(3)$-valued function~\eqref{eq:so3-target}, and the samples are corrupted by tangent-space Gaussian noise with standard deviation $\sigma$. Solid lines denote the median over the evaluation grid and shaded regions show the 25th--75th percentile over the evaluation grid, computed from the pointwise quantities estimated over 100 trials.}
\label{fig:so3}
\end{figure}

Based on the manifold extension of the BVT in~\eqref{eq:mse-manifold}--\eqref{eq:var-manifold}, we compare the Shepard approximation with its multiscale variant. The results are shown in Figure~\ref{fig:so3}. In terms of bias ratio, the multiscale approach consistently outperforms the standard Shepard approximation across all noise levels. Furthermore, Figure~\ref{fig:so3-mse} demonstrates that the multiscale method achieves a superior error rate, even under relatively high noise levels.

\subsection{Approximating a function with values on the SPD manifold}\label{sec:spd}

We further extend our previous experiment into a symmetric positive definite (SPD) manifold setting, which is a special case of a Hadamard manifold~\cite{bhatia2009positive}. We define the multiscale operators on the manifold as in~\cite {sharon2022multiscale} and use the affine-invariant Riemannian metric, with the corresponding closed-form exponential, logarithm, and geodesic distance. The noiseless target is defined in the tangent space of symmetric matrices and then projected to the SPD manifold:
\begin{equation}\label{eq:spd-target}
    G_0(x,y) = \exp_\M\left(\Psi(x,y)\right),
\end{equation}
where
\begin{equation*}
    \Psi(x,y)= \begin{pmatrix}
    \sin(2\pi y)\cos(2\pi x) & y^2 & xy\\
    y^2 &  1 & 0 \\ xy & 0 & \cos(\pi x)
    \end{pmatrix}.
\end{equation*}
Noisy samples are generated in the same tangent coordinates by
\[
    G_i^\varepsilon
    =
    \exp_\M\left(\Psi(x_i,y_i)\odot (J_{3\times3}+\sigma Z_i)\right).
\]
Here, $\exp_\M$ is the matrix exponential, $\odot$ is the Hadamard element-wise product, and $J_{3\times3}=[1]_{3\times3}$. Each~$Z_i$ is a symmetric random matrix generated by sampling $(Z_i)_{ab}\sim\mathcal N(0,1)$ independently for $a\leq b$ and setting $(Z_i)_{ba}=(Z_i)_{ab}$; the symbol~$Z$ avoids a clash with the noise covariance~$\Sigma$ of Assumption~\ref{ass:noise}. The parameter~$\sigma$ is a fixed noise amplitude. As an expression for the SNR, we use the empirical mean of~$\operatorname{SNR}(x_i,y_i)={\norm{\Psi(x_i,y_i)}_F^2}/{\norm{\Psi(x_i,y_i)\odot \sigma Z_i}_F^2}$, measured with the Frobenius norm~$\norm{\cdot}_F$. Thus, the reported SNR quantifies the perturbation size in the tangent space. The dataset consists of about $11^2$ points, sampled uniformly at random from the square $[0,1]^2$.

\begin{figure}[!htb]
\centering
\begin{subfigure}{.45\textwidth}
\centering
\includegraphics[width=\textwidth]{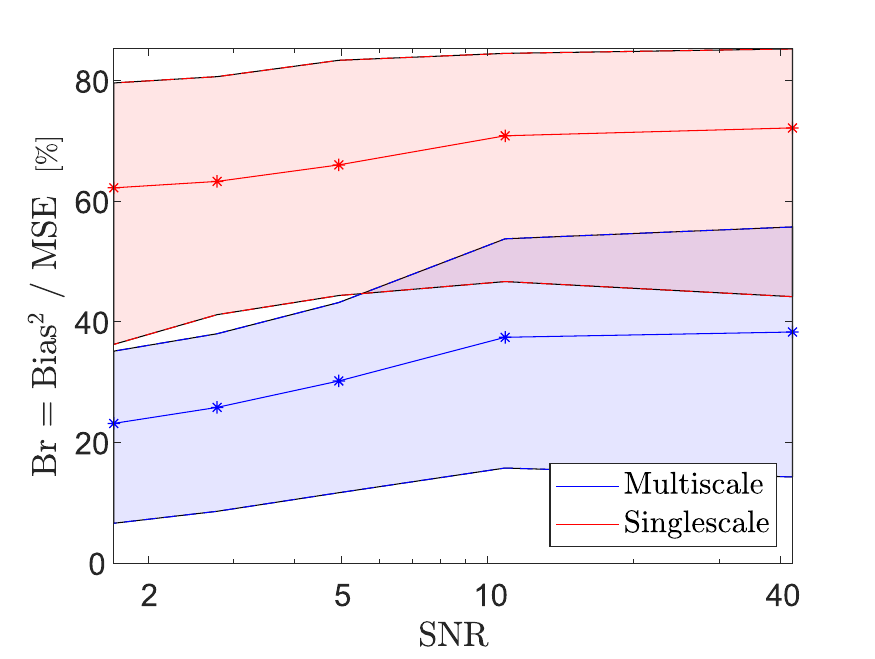} 
\centering
\caption{Bias ratio}\label{fig:spd-br}
\end{subfigure} \quad
\begin{subfigure}{.418\textwidth}
\includegraphics[width=\textwidth]{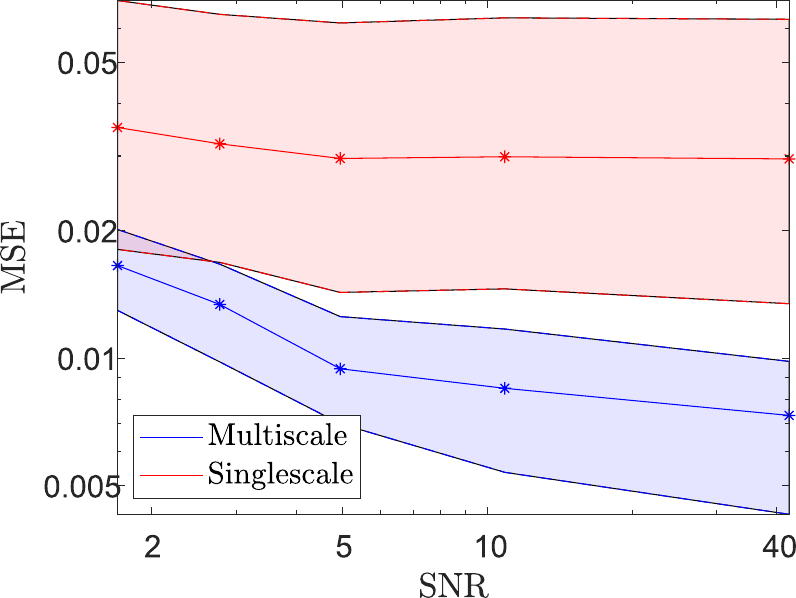}
     \centering   
    \caption{Error rates}
    \label{fig:spd-mse}
\end{subfigure}
\caption{Comparing Shepard (red) and its multiscale (blue) across different signal-to-noise ratios. The target function is the noiseless SPD-valued function~\eqref{eq:spd-target}, and the samples are corrupted by tangent-space noise corresponding to the reported SNR. Solid lines denote the median over the evaluation grid and shaded regions show the 25th--75th percentile over the evaluation grid, computed from the pointwise quantities estimated over 100 trials.}
\label{fig:spd}
\end{figure}

As in the previous examples, Figure~\ref{fig:spd-br} shows that the multiscale scheme consistently yields a lower bias ratio. In Figure~\ref{fig:spd-mse}, the multiscale approach also yields superior error rates.

\subsection{Diffusion operators}\label{sec:diffusion}

This experiment directly illustrates the spectral model of
Theorem~\ref{thm:spectral-damping}. We fix an equispaced grid
$X=\{x_1,\ldots,x_N\}\subset[0,1]^2$ with spacing $0.075$
and $N=196$ sites. All approximations and error measurements
are evaluated on this fixed grid. The integrated quantities
are computed using the Euclidean norm on $\mathbb R^N$,
that is, by summing over the data sites, with no normalization by $N$.

We construct a symmetric graph Laplacian $L=D_W-W$, where
\[
W_{ab}=
\begin{cases}
\exp\!\left(-\dfrac{\|x_a-x_b\|_2^2}{2\ell^2}\right),
    & a\ne b,\\
0,  & a=b,
\end{cases}
\qquad
D_W=\operatorname{diag}(W\mathbf1),
\]
and $\ell$ is the median of the pairwise distances.
Writing $L=U\Lambda U^\top$, we define three diffusion filters
\[
Q_i=U\exp(-\tau_i\Lambda)U^\top,
\qquad
\tau_i=\frac{1}{10\cdot2^i},
\qquad i=1,2,3.
\]
Each level acts on the full data vector, with the scale determined by the diffusion time $\tau_i$. Thus, all operators share the orthonormal eigenbasis of $L$ and have spectral responses in $[0,1]$, as required by the theorem. The multiscale estimator is computed as
\[
\widetilde f_j
=
\left[I-(I-Q_j)\cdots(I-Q_1)\right]Y,
\qquad j=1,2,3,
\]
and is compared with the single-scale estimator $Q_jY$.

For a fixed ordering of the sites and of the eigenvectors, the noiseless target vector is prescribed through its spectral coefficients,
\[
f=U\theta,
\qquad
\theta_k=\|x_k\|_2,
\qquad k=1,\ldots,N,
\]
The random site ordering assigns the coefficients without imposing decay with increasing eigenvalue. The target therefore tests the spectral formulas with energy distributed throughout the spectrum. The damping of each mode is determined separately by the filter responses $\lambda_{i,k}=\exp(-\tau_i\omega_k)$. The analytical curves below use these same target coefficients and spectral responses.

In each of $T=100$ independent trials, we generate
\[
Y^{(t)}=f+\sigma z^{(t)},
\qquad z^{(t)}\sim\mathcal N(0,I_N).
\]
The sites, target, and diffusion operators remain fixed across trials. Consequently, the noise has zero mean and covariance $\sigma^2I_N$, exactly as assumed in Theorem~\ref{thm:spectral-damping}. As in the preceding plots, we use $\operatorname{SNR}=1/\sigma$ as a noise-scale parameter, with the noiseless case marked separately.

We estimate MISE by averaging the squared Euclidean errors
over $T=100$ independent trials, and estimate the integrated
squared bias by $\|\frac{1}{T}\sum_{t=1}^{T}\widetilde f_j^{(t)}-f\|_2^2$. Dividing the estimated integrated squared bias by the estimated MISE gives the empirical integrated bias ratio. The analytical curves are computed from the exact formulas
of Theorem~\ref{thm:spectral-damping}, using the same target
coefficients, spectral responses, and noise variance.
The empirical quantities are finite-trial estimates of these
population quantities; in particular, the expected squared-bias
estimate equals $\ibias_j+\ivar_j/T$.

\begin{figure}[!ht]
\centering
\begin{subfigure}{.47\textwidth}
    \centering
    \includegraphics[width=\textwidth]{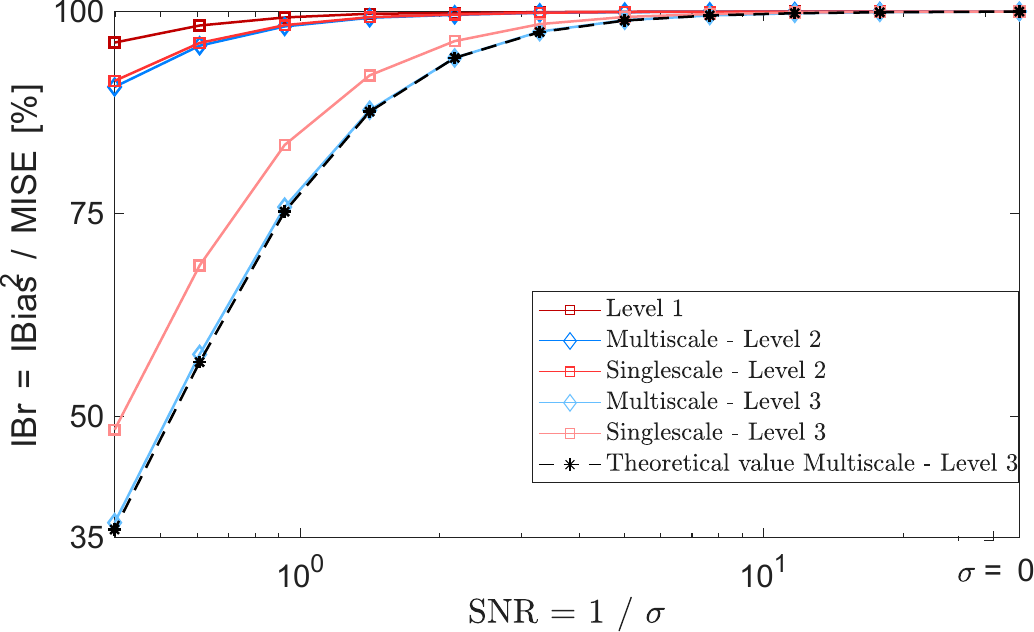}
    \caption{Integrated bias ratio}
    \label{fig:diffusion-br}
\end{subfigure}
\quad
\begin{subfigure}{.47\textwidth}
    \centering
    \includegraphics[width=\textwidth]{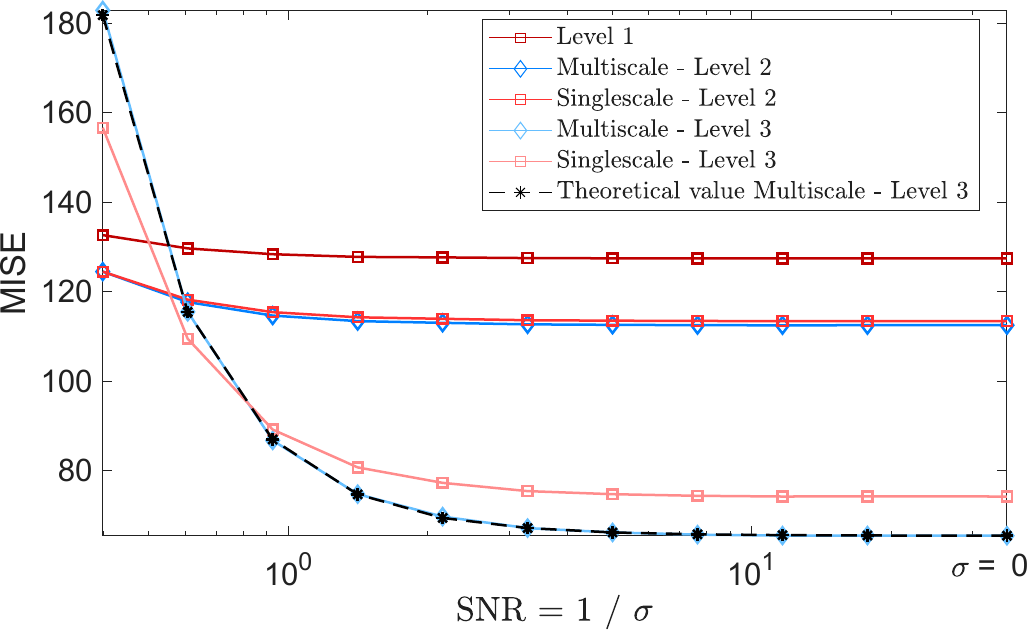}
    \caption{MISE}
    \label{fig:diffusion-mise}
\end{subfigure}
\caption{Integrated bias ratio (left) and MISE (right) for
single-scale diffusion estimators (red) and their multiscale
counterparts (blue), across noise levels. The experiment uses
fixed full-data diffusion filters and additive isotropic
Gaussian noise, satisfying the assumptions of
Theorem~\ref{thm:spectral-damping}. Empirical quantities are
estimated from 100 independent trials. Black curves show
the exact analytical values for the third multiscale level.}
\label{fig:diffusion}
\end{figure}

Figure~\ref{fig:diffusion} compares the single-scale and multiscale estimators with the analytical predictions, shown for the
third multiscale level. The close agreement supports the
predicted bias--variance behavior. Additional levels reduce
the integrated bias ratio, while their effect on MISE depends
on the noise level: the variance increase can outweigh the
bias reduction. This illustrates the distinction between
bias-ratio monotonicity and error reduction quantified by
Corollary~\ref{cor:mise-increment}.

\section{Conclusion}

This study demonstrates the effectiveness of multiscale approximation methods in significantly reducing bias while improving overall approximation accuracy. The theoretical results in Section~\ref{sec:residual-bias} identify the mechanism behind this behavior: in the scalar-valued linear setting, the bias of the multiscale estimator is exactly the deterministic residual~$\mathcal R_nf$, and a residual-contraction condition therefore implies bias reduction. Under a spectral-filter model, the residual is damped componentwise, and the integrated bias ratio is monotone nonincreasing with the number of multiscale levels; moreover, an exact per-level error increment exhibits a bias--variance crossover that provides a criterion for when an additional level is beneficial; the resulting number of levels depends on the spectral responses of the level operators and on the signal-to-noise ratio, and it is not universally small. A local tangent-space argument further explains why the same residual-bias mechanism persists to leading order for manifold-valued quasi-interpolation based on intrinsic averages, and yields a leading-order bias--variance decomposition of the manifold-valued error.

Through numerical experiments, we confirmed that this mechanism is visible beyond the idealized assumptions of the theory. The experiments demonstrate reductions in bias and, in many of the tested settings, improvements in MSE; the diffusion experiment also illustrates when the variance increase outweighs the bias reduction. The bias ratio~\eqref{eq:bias-ratio} provides a straightforward, interpretable metric for comparing the bias content of different operators and settings at a fixed number of levels, provided it is read together with the MSE. As Remark~\ref{rem:optimal-levels} makes explicit, it is deliberately not a criterion for selecting the number of levels, which remains an MSE consideration governed by the crossover criterion~\eqref{eq:crossover}.

Future research could combine the deterministic convergence rates recently obtained for multilevel MLS on regular grids~\cite{durst2026improved} with the residual--bias identity of Theorem~\ref{thm:bias-equals-residual} to obtain quantitative bias rates, and extend such rates from gridded to scattered data sites; further directions are to sharpen the curvature-dependent bounds in the manifold-valued setting, extend the analysis to broader classes of nonlinear quasi-interpolants, and investigate the impact of reduced bias on advanced applications, such as data representation and geometric learning~\cite{rabin2012heterogeneous}, and model order reduction~\cite{zimmermann2021manifold, farhat2020computational}. Overall, our findings suggest that multiscale quasi-interpolation offers a simple yet powerful tool for bias control in numerical modeling, with potential applications in data-driven geometric analysis and manifold learning.

\section*{Acknowledgements}
NS is partially supported by the NSF-BSF award 2024791, the BSF award 2024266, and the DFG award 514588180.

\section*{CRediT authorship contribution statement}
\textbf{Asaf Abas:} Conceptualization, Methodology, Software, Validation, Investigation, Visualization, Writing -- original draft. \textbf{Nir Sharon:} Conceptualization, Methodology, Supervision, Funding acquisition, Writing -- review \& editing.

\section*{Declaration of competing interest}
The authors declare no known competing financial interests or personal relationships that could have influenced the work reported in this paper.

\section*{Declaration of generative AI and AI-assisted technologies}
During the preparation of this work, the authors used ChatGPT (OpenAI) in order to assist with manuscript organization and language revision. After using this service, the authors reviewed and edited the content as needed and take full responsibility for the content of the publication.

\section*{Data availability}
The code that reproduces the numerical experiments is publicly available at \url{https://github.com/ABASASA/Multiscale_Bias_Ratio}.

\bibliographystyle{plain} 
\bibliography{Multiscale.bib}

\end{document}